\documentclass{article}
\usepackage[utf8]{inputenc}
\usepackage{amsmath}
\usepackage{amssymb}
\usepackage{amsthm}
\DeclareMathOperator*{\argmin}{arg\,min}
\usepackage{array}
\usepackage{float}
\usepackage{geometry}
\usepackage{pgfplots}
\usepackage[authoryear]{natbib}
\usepackage{dsfont}
\usepackage[toc,page]{appendix}
\usepackage{mathrsfs}
\usepackage{hyperref}
\usepackage[ruled,linesnumbered]{algorithm2e}
\usepackage{etoolbox}
\AtBeginEnvironment{algorithm}{\SetInd{.5em}{1em}}
\SetKwInput{KwInput}{Input}                
\SetKwInput{KwOutput}{Output}              

\newtheorem{proposition}{Proposition}
\newtheorem{definition}{Definition}

\newcommand{\probP}{\text{I\kern-0.15em P}}
\newcommand{\vrpsd}{VRPSD\xspace}

\usepackage{authblk}

\title{The disaggregated integer L-shaped method for the stochastic vehicle routing problem}

\author[(1)]{Lucas Parada}
\author[(2)]{Robin Legault}
\author[(3)]{Jean-Fran\c{c}ois C\^ot{\'e}\thanks{Corresponding author}}
\author[(1)]{Michel Gendreau}

\affil[(1)]{CIRRELT, Polytechnique Montr\'eal,} 
\affil[  ]{\tt lucas.parada@polymtl.ca}
\affil[  ]{\tt michel.gendreau@polymtl.ca}
\affil[(2)]{CIRRELT, Université de Montr\'eal,}
\affil[  ]{\tt robin.legault@umontreal.ca}
\affil[(3)]{CIRRELT, Universit\'e Laval,}
\affil[  ]{\tt jean-francois.cote@fsa.ulaval.ca}
\date{}

\begin{document}

\maketitle
\begin{abstract}
This paper proposes a new integer L-shaped method for solving two-stage stochastic integer programs whose first-stage solutions can decompose into disjoint components, each one having a monotonic recourse function. In a minimization problem, the monotonicity property stipulates that the recourse cost of a component must always be higher or equal to that of any of its subcomponents. The method exploits new types of optimality cuts and lower bounding functionals that are valid under this property. The stochastic vehicle routing problem is particularly well suited to be solved by this approach, as its solutions can be decomposed into a set of routes. We consider the variant with stochastic demands in which the recourse policy consists of performing a return trip to the depot whenever a vehicle does not have sufficient capacity to accommodate a newly realized customer demand. This work shows that this policy can lead to a non-monotonic recourse function, but that the monotonicity holds when the customer demands are modeled by several commonly used families of probability distributions. We also present new problem-specific lower bounds on the recourse that strengthen the lower bounding functionals and significantly speed up the resolution process. Computational experiments on instances from the literature show that the new approach achieves state-of-the-art results.

\end{abstract}

\section{Introduction}

In two-stage stochastic programs, first-stage decisions are taken before the realization of uncertain parameters. Once these parameters are revealed, second-stage decisions are taken to correct the first-stage solution. The objective is to find a first-stage solution that minimizes the cost of the first stage and the expected cost of the second stage. These problems have typically been solved using Benders decomposition \citep{benders1962partitioning}, or equivalently the L-shaped method \citep{VanSlyke1969} when they have continuous variables and a large number of scenarios. For solving programs with binary variables in both stages, researchers historically relied on the integer L-shaped method of \cite{laporte1993integer}. These methods can be regarded as branch-and-cut (B\&C) algorithms where the expectation function of the objective is relaxed and replaced by a variable that is linearly bounded with inequalities. Nowadays, the L-shaped and integer L-shaped methods are considered to be computationally impractical because of their slow convergence. Several techniques have been proposed in the literature to mitigate this problem. A common approach is to add lower bounding functionals (LBFs) \citep{Birge1986} to the model. These inequalities bound the recourse for a much broader range of solutions than optimality cuts. However, many of the most impactful ones are problem-specific, which limits their range of application. Other methods like the logic-based Benders decomposition of \cite{hooker2003logic} try to overcome the problem by providing cuts that include only a small subset of variables. These cuts are then effective for a large number of solutions. Other enhancements are reported in \cite{RAHMANIANI2017801}.

This paper introduces a new integer L-shaped method, named the disaggregated integer L-shaped method (DL-shaped method). This method is tailored for a specific class of two-stage stochastic programs in which, given a first-stage solution, the recourse cost can be expressed as a sum of independent recourse functions involving disjoint sets of first-stage variables. The DL-shaped method also requires each of the resulting recourse functions to be monotonic. In the case of a minimization problem, this property means that adding first-stage variables to a component cannot decrease its recourse cost. Many two-stage stochastic programs have such structure and property. For example, in the bin packing problem with uncertain weights, it is possible to decompose the first-stage solution by bin. This problem is found in the allocation of surgeries to operating rooms \citep{denton2010}. Similarly, in job scheduling problems \citep{Li2022, elci2022}, jobs are assigned to machines, and their first-stage solution can be decomposed by machine. In the latter application, the monotonicity of the objective function corresponds to the fact that assigning more jobs to a machine cannot decrease the expected makespan of its operations.

We propose an implementation of the DL-shaped method for the vehicle routing problem with stochastic demands (\vrpsd). In this problem, each customer has to be visited exactly once by a fleet of capacitated vehicles. This structure makes it possible to decompose the recourse cost by route. Customer demands are modeled by independent random variables (RVs) with known distributions and are only observed when a vehicle arrives at the customer's location. A failure occurs when a vehicle arrives at a customer's location with a remaining capacity smaller than the customer's demand. When this happens, a recourse action must be implemented to satisfy the current and next customer demands. This work applies the detour-to-depot (DTD) recourse action, which performs a back-and-forth trip from the customer to the depot to restock the vehicle. We will show that this recourse action has the property of monotonicity under some conditions.

The contributions of this paper are as follows. First, we introduce the DL-shaped method in the context of the \vrpsd. Second, we demonstrate that, although the DTD recourse function is not monotonic in general, it is for Poisson distributions and some normal, binomial, Erlang, and negative binomial distributions when the sum of the expected demands on each route respects the vehicle capacity. Third, we present several new lower bounds on the DTD recourse function that take advantage of its monotonicity and are used by the LBFs. To our knowledge, we perform the first numerical comparison of the different lower bounds on the recourse. Our numerical experiments show our new lower bounds greatly improve over the ones found in the literature and that our new method achieves state-of-the-art results on existing benchmark instances. This leads us to propose new instances that are more challenging than the existing ones from the literature.

The remainder of this paper is organized as follows. Section \ref{sec:Related literature} presents the literature on stochastic vehicle routing problems. Section \ref{sec:Mathematical formulations} presents the mathematical formulation of the problem. Section \ref{sec:Monotonicity} presents theoretical results on the monotonicity of the recourse function. Section \ref{sec:Disaggregated integer L-shaped method} details the implementation of the new method as well as the new lower bounds for the problem. Section \ref{sec:Numerical results} presents the numerical results for our implementation, and lastly, Section \ref{sec:Conclusions} presents the conclusion.

\section{Related literature} \label{sec:Related literature}

The stochastic vehicle routing problem with random customer demands and multiple depots by \cite{Tillman1969} is considered to be the first studied version of the \vrpsd. Nowadays, the state-of-the-art exact methods can solve instances with more than 100 customers \citep{florio2020new}. By this measure, since \cite{Tillman1969}, the field has linearly increased over time in its solution capability. Although this is considerable progress for an NP-Hard problem like the \vrpsd, \cite{gendreau201650th} conclude that we are still in the early stages of the field, with room for significant improvement in exact solution methods. In this section, we present a historical review of the exact methods developed for the \vrpsd.

The a priori and re-optimization paradigms are the two modeling paradigms used in the literature. The choice of either of these paradigms depends on when the stochastic parameters are revealed in the solution process. Formalized by \cite{bertsimas1990priori}, the a priori paradigm relies on the design of routes that cannot be modified after they have been fixed. Re-optimization was formalized by \cite{dror1989vehicle} as a Markov decision process (MDP) and allows for the construction and adjustment of routes as information is gradually revealed. In recent decades, the increasing availability of information technologies has simplified the gathering of real-time information that helps decision-makers re-optimize a priori routes on the fly. As a consequence, new paradigms have emerged, such as the dynamic and online vehicle routing problems, where the parameter are updated throughout the solution process. A comprehensive survey on the dynamic vehicle routing problem is presented in \cite{PILLAC20131}.

The DTD and optimal restocking (OR) policies are the most common in the literature. Under the DTD policy, restocking trips are only allowed in case of failure. Studies of the \vrpsd with this policy include the works of \cite{gendreau1995exact}, \cite{laporte2002integer}, and \cite{jabali2014partial}. For the OR policy, vehicles are allowed to perform preventive returns to the depot to avoid costly failures. \cite{yangRestocking} have shown that the OR policy can be calculated by solving a Bellman equation. Relevant works under this policy include \cite{louveaux2018exact}, \cite{salavati2019exact}, \cite{florio2020new}, and \cite{Ymro2023}. Additionally, variants of these two policies have also been studied. For example, \cite{Aykagan2007} propose a recourse policy that is based on pairwise sets of a priori routes. Also, in \cite{Secomandi2009Reoptimization}, an a priori planned tour is partly re-optimized using a MDP that is limited to a restricted set of states. \cite{florio2022vehicle} apply another type of partial re-optimization in which the order of visit of successive customers in a planned route is allowed to be switched.

Several numerical experiments have been performed in the literature to compare the effectiveness of the proposed policies. Those of \cite{louveaux2018exact} show that the OR policy might yield some savings compared to the DTD policy, but they also noticed that the optimal a priori routes are most of the time the same as those from the DTD policy. Their results also indicate that more instances can be solved when using the DTD policy. \cite{salavati2019exact} computed the cost of the DTD policy from an optimal OR policy solution. An average increase in the optimal value of about 0.63\% could be observed. \cite{florio2020new} indicate that the benefit of using the OR policy over the DTD policy is close to none. However, savings of 1.6\% on average can be obtained when reducing the number of possible routes and allowing the capacity to be exceeded. They also noted that the benefit of using the OR policy decreases when the demand variability increases. Finally, the partial re-optimization policy of \cite{florio2022vehicle} can further reduce costs compared to the OR policy by 0.73\% on average. There are still open questions from these observations. Those analyses are limited to relatively few instances of small sizes. In light of these results, the DTD policy is still worth studying.

Exact methods from the literature for the \vrpsd divide into two categories, namely integer L-shaped and branch-and-price (B\&P) methods. For the integer L-shaped method, its main advantages are that it can manage solutions that include long routes (40+ nodes per route), it is compatible with several demand distributions, including continuous ones, and it can handle many different types of recourse functions. Unfortunately, it is less efficient for instances requiring more than three vehicles. Although several authors have devoted important efforts to propose LBFs \citep{laporte2002integer,jabali2014partial,Ymro2023} that alleviate this limitation, solving instances with a large number of vehicles remains a challenge. On the other hand, B\&P methods have shown to be efficient for instances that require a large number of short routes, but they also have important limitations. They have been almost exclusively implemented for demands following Poisson distributions. Also, for the OR policy, they require to store entire routes in the labeling algorithm of the column generation subproblem, which greatly limits the length of the routes the method can handle. Finally, for the DTD policy, they do not scale well when the demands and the capacity of vehicles are high since the size of the pricing subproblem depends on the vehicle capacity. To mitigate this limitation, some methods from the literature divide each instance's capacity and expected demands by their greatest common divisor to produce subproblems of tractable size in their column generation procedures. Although the instances resulting from this transformation are usually good approximations of the original problems, they are not equivalent in general and can have different optimal solutions and values.

The idea of disaggregating the recourse, i.e., separating it into independent disjoint functions, was studied by \cite{wets1966programming}. They show that the recourse can be disaggregated trivially when the second-stage linear programming matrix can be decomposed into blocks. In such cases, a variable is added in the first stage for each block. Then, in each iteration, the second-stage problem of each block is solved, and an optimality cut that contains only first-stage variables of the block is added. Extending this idea to problems whose second-stage matrix is not decomposable is challenging.

The second stage of vehicle routing problems does not present this block structure. One of the reasons for this is that the uncertain parameters are progressively revealed as the clients are visited. However, the first-stage solutions of vehicle routing problems have a separable structure, as they form independent routes. Some works have introduced techniques to exploit this separability. \cite{seguin1996problemes} proposed to disaggregate the recourse by routes, \cite{cote2020vehicle} proposed to bound the recourse for sets of nodes, and \cite{Ymro2023} for unstructured sets. These methods add variables in the first stage to bound some parts of the second stage by using linear inequalities. Numerical experiments show that these valid inequalities help at reducing computation times. The superiority of the DL-shaped method over these previous works comes from its new optimality cuts and new LBFs that simultaneously bound the recourse for paths and sets of nodes that apply to a broad range of first-stage solutions.

As noted in the introduction, the DL-shaped method requires the monotonicity of the recourse function. The topic of monotonicity received little attention in the \vrpsd literature. To the best of our knowledge, the articles by \cite{dror1987InventoryRouting} and \cite{KREIMER199063} constitute the only theoretical works on the matter. They both give sufficient conditions for the probability of a failure to increase at each customer on a route. Although we show that this property is sufficient for the recourse function to be monotonic under the DTD policy, the results of \cite{dror1987InventoryRouting} and \cite{KREIMER199063} only apply to independent and identically distributed (i.i.d.) demands. This paper introduces sufficient conditions for the monotonicity of the recourse function that can be applied in the independent but not identically distributed case. We show conditions under which it is respected for Poisson, normal, binomial, Erlang, and negative binomial demands.

\section{Mathematical formulation} \label{sec:Mathematical formulations}

This section introduces the notation that will be used throughout the paper and the mathematical model of the problem.

The \vrpsd is defined as follows. Let $G=\{N_0,E\}$ be an undirected graph where $N_0=\{0,1,...,n\}$ is the set of nodes, with node $0$ being the depot, $N=\{1,...,n\}$ the set of customer nodes, and $E=\{(i,j):i,j \in N, i < j\}$ the set of edges. Traveling on edge $(i,j)$ incurs a travel cost of $c_{ij}$. Each customer $i\in N$ has a non-negative demand given by the RV $\xi_i$, with an expected value $\mu_i$ and a standard deviation $\sigma_i$. It is assumed that the demand variables are independently distributed. A fleet of identical vehicles is available to satisfy customer requests. Each vehicle has a capacity $Q$ and must follow a route that starts and ends at the depot. Each customer must be visited exactly once. A route is defined as a sequence of customers starting and ending at the depot, and a sequence is feasible if the sum of the expected demands respects the vehicle capacity. This condition was added by \cite{laporte2002integer} to avoid routes that systematically fail while having others that are underutilized. Let $M \subseteq \mathds{N}$ be the set specifying the number of vehicles that can be used in a solution. For unlimited fleet problems, this set is given by $M=\{ \lceil (\sum_{i \in S}\mu_i)/Q \rceil, \dots, |N|\}$, and when exactly $\bar m$ vehicles must be used, $M=\{\bar{m}\}$.

Let $E(S)$ be the edges with both endpoints in $S$ and $E(h)$ the set of edges linked with node $h$. The \vrpsd can be formulated as follows. The variable $x_{ij}$ indicates the number of times edge $(i, j) \in E$ is traversed. The variables associated with pairs of customers are binary, and those linking the depot to a customer can take the values 0, 1, or 2. The variable $x_{0i}$ equals 2 if a vehicle serves a single customer $i \in E$. A binary variable $z_m$ decides whether $m$ vehicles are used in the solution. Let $\mathcal{Q}(x)$ denote the expected recourse cost of the first-stage solution $x = (x_{ij})$. The model is:
\begingroup
\allowdisplaybreaks
\begin{alignat}{2}
    \min & \sum_{(i,j)\in E}c_{ij}x_{ij} + \mathcal{Q}(x) \label{model:obj} \\
    \text{s.t. } & \sum_{i \in N} x_{0i} = \sum_{m \in M} 2mz_m, \label{model:depot_degree} \\
    & \sum_{(i,j) \in E(h)} x_{ij} = 2 & h \in N, \label{model:degree}\\
    & \sum_{(i,j) \in E(S)} x_{ij} \leq |S| - \left\lceil \frac{\sum_{i \in S}\mu_i}{Q} \right\rceil  & S \subseteq N, \label{model:cap}\\
    & \sum_{m \in M} z_m = 1, \label{model:vehicle_eq_1} \\
    & x_{ij} \in \{0,1\} &\qquad (i,j) \in E(N) \label{model:bin},\\
    & x_{ij} \in \{0,1,2\} & (i,j) \in E(0). \label{model:int}
\end{alignat}
\endgroup

The objective \eqref{model:obj} is to minimize the sum of travel costs plus the expected recourse cost. Constraint \eqref{model:depot_degree} imposes that $m$ routes must be connected to the depot if $m$ vehicles are used. Constraints \eqref{model:degree} ensure that customers are visited exactly once. Constraints \eqref{model:cap} impose that the expected demand on each route does not exceed the vehicle capacity and that the routes are connected to the depot. Constraint \eqref{model:vehicle_eq_1} ensures that the number of vehicles that are used in the solution is an element of $M$. Constraints (\ref{model:bin}) and \eqref{model:int} define the domain of the variables.

For the DTD recourse policy, it is possible to separate the calculation of the recourse by route. Let $\mathcal{R}^{\nu}$ be the set of routes in solution $x^\nu$. Then, the expected recourse cost is defined as:
\begin{align}
    & \mathcal{Q}(x^\nu)=\sum_{r \in \mathcal{R}^{\nu}} \mathcal{Q}(r),\label{eq:8}
\end{align}
where $\mathcal{Q}(r)=\min\{\mathcal{Q}^1(r),\mathcal{Q}^2(r)\}$ is the minimum expected recourse cost of both orientations of route $r = (0,i_1,\ldots,i_t,0)$. The computation of both orientations is required because model (\ref{model:obj})-(\ref{model:int}) does not capture routes orientation. The expected recourse cost is denoted by $\mathcal{Q}^1(r)$ and $\mathcal{Q}^2(r)$ in the two orientations. Under the DTD strategy, a cost of $2c_{0i}$ is incurred when the demand of customer $i$ cannot be satisfied. The expected recourse of route $r$, in the first orientation, is thus calculated as follows.
\begin{equation}
  \label{eq:rec_route}
  \mathcal{Q}^1(r)= 2\sum_{j=1}^t\sum_{l=1}^{\infty}\mathds{P} \left ( \sum_{k=1}^{j-1} \xi_k \leq lQ < \sum_{k=1}^{j} \xi_k \right)c_{0j}
\end{equation}

Note that a restocking trip occurs only when a customer's demand exceeds the vehicle's residual capacity. We exclude the case of going to the depot and then moving to the next customer when an exact stock-out occurs. If it is assumed that $\xi_i \leq Q$ for all $i$ with probability 1, then the second summation of (\ref{eq:rec_route}) can be limited to $j-1$ instead of going to $\infty$. If demands can exceed the capacity of the vehicles, the probability $\mathds{P} \left( \sum_{k=1}^{j-1} \xi_k \leq lQ \leq \sum_{k=1}^{j} \xi_k \right)$ tends to zero as $l$ increases and the summation can be stopped when it reaches a probability close to zero. In our implementation, the summation is stopped when the probability is lower than 0.0001. Throughout the paper, the expected recourse cost $\mathcal{Q}(r)$ of a route $r=(0,p,0)$ will equivalently be denoted by $\mathcal{Q}(p)$.

\section{Monotonicity of the recourse function}\label{sec:Monotonicity}

Our algorithm relies on new optimality cuts and lower bounds that require the recourse function to decrease monotonically when removing customers from a route. This section formalizes this property by the so-called monotonicity condition and monotonicity property, which respectively apply to sets of customers and instances of the \vrpsd. Based on these definitions, we prove sufficient conditions for the recourse function to be monotonic under the DTD policy when the demands follow independent Poisson, normal, binomial, Erlang, and negative binomial distributions. We also provide counterexamples to illustrate that the recourse function is not monotonic in general.

\begin{definition}
    Let $p=(i_1, \dots, i_t)$ be a path on graph $G$. We say that $p'=(j_1, \dots, j_{t'})$ is a subpath of $p$ if it can be obtained by removing elements of $p$, without changing the ordering of the remaining elements. Formally, $p'$ is thus said to be a subpath of $p$ if $\{j_1, \dots, j_{t'}\} \subseteq \{i_1, \dots, i_t\}$ and, for any pair $(j_{a'}, j_{b'}) = (i_a, i_b)$, $a' < b'$ if and only if $a < b$.
\end{definition}

\begin{definition}
    Let $p=(i_1, \dots, i_t)$ be a path on graph $G$. A subpath $p'$ that can be written as $p'=(i_a, \dots, i_b)$ for some $1\leq a \leq b \leq t$ is said to be connected. In other words, a connected subpath $p'$ is obtained by removing elements at the beginning and the end of $p$.
\end{definition}

\begin{definition}\label{def:monotonicity_condition}
    A set of customers $S \subseteq N$ is said to respect the monotonicity condition if, for any pair $(a, b) \in S \times S$ such that $a \neq b$ and for any subset $\tilde{S} \subseteq S \setminus \{a, b\}$, the following inequality:
    \begin{equation}
        \label{ineq:monotonicity_condition_1} \mathds{P}\left(\tilde{\xi} + \xi_a \leq lQ <\tilde{\xi} + \xi_a + \xi_b \right) \geq \mathds{P}\left(\tilde{\xi} \leq lQ <\tilde{\xi} + \xi_b\right),
    \end{equation}
where $\tilde{\xi} = \sum_{i \in \tilde{S}} \xi_{i}$, holds for any positive integer $l \in \mathds{N}$.
\end{definition}
\begin{definition}\label{def:monotonicity_property}
    An instance of the \vrpsd with the DTD policy is said to have the monotonicity property if any set $S \subseteq N$ such that $\sum_{i \in S} \mu_i \leq Q$ respects the monotonicity condition.
\end{definition}

The above definition can be applied by enumeration to verify whether an instance respects the monotonicity property, but it might be unpractical for large instances since the number of comparisons grows exponentially with $|N|$. In Propositions \ref{prop:monotonicity_condition_2} and \ref{prop:monotonicity_condition_sufficient}, we introduce an equivalent formulation of the monotonicity condition and a sufficient condition which will later be used to prove that the monotonicity property is respected if the demands are modeled by specific distributions.

\begin{proposition}
\label{prop:monotonicity_condition_2}
Inequality (\ref{ineq:monotonicity_condition_1}) can be rewritten equivalently as follows.
\begin{equation}
    \label{ineq:monotonicity_condition_2} \mathds{P}\left(\tilde{\xi} + \xi_a \leq lQ\right) - \mathds{P}\left(\tilde{\xi} + \xi_a + \xi_b \leq lQ\right) \geq \mathds{P}\left(\tilde{\xi} \leq lQ\right) - \mathds{P}\left(\tilde{\xi} + \xi_b \leq lQ\right)
\end{equation}
\end{proposition}
\proof{}
\begin{align}
    \label{monotonicity_condition_e1} \mathds{P}\left(\tilde{\xi} + \xi_a \leq lQ <\tilde{\xi} + \xi_a + \xi_b \right) &\geq \mathds{P}\left(\tilde{\xi} \leq lQ <\tilde{\xi} + \xi_b\right)\\
    \label{monotonicity_condition_e2} \iff \mathds{P}\left(\tilde{\xi} + \xi_a + \xi_b > lQ \right) -  \mathds{P}\left(\tilde{\xi} + \xi_a > lQ \right) &\geq \mathds{P}\left(\tilde{\xi} + \xi_b > lQ \right) -  \mathds{P}\left(\tilde{\xi} > lQ \right)\\
    \iff \mathds{P}\left(\tilde{\xi} + \xi_a \leq lQ \right) - \mathds{P}\left(\tilde{\xi} + \xi_a + \xi_b \leq lQ \right) &\geq \mathds{P}\left(\tilde{\xi} \leq lQ \right) -\mathds{P}\left(\tilde{\xi} + \xi_b \leq lQ \right)
\end{align}
Inequalities (\ref{monotonicity_condition_e1}) and (\ref{monotonicity_condition_e2}) are equivalent since it is assumed that the demands have non-negative support.
\endproof

\begin{proposition}
\label{prop:monotonicity_condition_sufficient}

Let $S \subseteq N$ be a set of customers such that, for any pair $(a, b) \in S \times S$, $a \neq b$, and for any subset $\tilde{S} \subseteq S \setminus \{a, b\}$, the random variables $\tilde{\xi}$, $\xi_a$ and $\xi_b$ can respectively be rewritten as the sum of $\tilde{\rho} \in \mathds{N}$, $\rho_a\in \mathds{N}$ and $\rho_b\in \mathds{N}$ i.i.d. non-negative RVs $\zeta_k$, $k \in \mathds{N}$. The following condition is sufficient for the set $S$ to respect the monotonicity condition:
\begin{equation}
    \label{ineq:monotonicity_condition_sufficient} P_{lQ}(\tilde{\rho}+\rho_a+j) \geq P_{lQ}(\tilde{\rho}+j) \hspace{0.5cm} \forall l \in \mathds{N}, \forall j \in \{1,\dots,\rho_b\},
\end{equation}
where $P_{u}(s) = \mathds{P}\left( \sum_{k=1}^{s-1} \zeta_k \leq u < \sum_{k=1}^{s} \zeta_k \right)$ denotes the probability that exactly $s \in \mathds{N}$ RVs $\zeta_k$ must be drawn for their cumulative sum to first exceed a given threshold $u \geq 0$.
\end{proposition}
\proof{}
Given the assumptions of the proposition, each original customer $i \in S$ can be seen as a group of $\rho_i$ smaller customers with i.i.d. demands $\zeta_k$, $k \in \mathds{N}$ and sharing the same location.

From there, the left-hand side of inequality (\ref{ineq:monotonicity_condition_1}) can be rewritten as follows.
\begin{align*}
    \mathds{P}\left(\tilde{\xi} + \xi_a \leq lQ <\tilde{\xi} + \xi_a + \xi_b\right) &= \mathds{P}\left( \sum_{k=1}^{\tilde{\rho}+\rho_a} \zeta_k \leq lQ < \sum_{k=1}^{\tilde{\rho}+\rho_a+\rho_b} \zeta_k \right) \\
    &= \sum_{j=1}^{\rho_b}\mathds{P}\left( \sum_{k=1}^{\tilde{\rho}+\rho_a+j-1} \zeta_k \leq lQ < \sum_{k=1}^{\tilde{\rho}+\rho_a+j} \zeta_k \right)\\
    &= \sum_{j=1}^{\rho_b} P_{lQ}(\tilde{\rho}+\rho_a+j)
\end{align*}

The second equality is obtained by decomposing the probability that the $l$-th restocking trip occurs between customers $\tilde{\rho}{+}\rho_a{+}1$ and $\tilde{\rho}{+}\rho_a{+}\rho_b$ as the sum of the probabilities that it occurs at customer $\tilde{\rho}{+}\rho_a{+}j$, for $j{=}1,\dots,\rho_b$.

Doing the same for the right-hand side, we obtain the sufficient condition for inequality (\ref{ineq:monotonicity_condition_1}).
\begin{align*}
    &\mathds{P}\left(\tilde{\xi} + \xi_a \leq lQ <\tilde{\xi} + \xi_a + \xi_b \right) \geq \mathds{P}\left(\tilde{\xi} \leq lQ <\tilde{\xi} + \xi_b\right)\\
    \iff &\sum_{j=1}^{\rho_b} P_{lQ}(\tilde{\rho}+\rho_a+j) \geq \sum_{j=1}^{\rho_b} P_{lQ}(\tilde{\rho}+j)\\
    \impliedby & P_{lQ}(\tilde{\rho}+\rho_a+j) \geq P_{lQ}(\tilde{\rho}+j) \hspace{2cm} \forall j \in \{1,\dots,\rho_b\} 
\end{align*}
\endproof

We now prove that, if the customers of a path respect the monotonicity condition, removing customers from this path cannot increase its recourse cost. This is shown in Proposition \ref{prop:monotonic_expected_recourse}.

\begin{proposition}
\label{prop:monotonic_expected_recourse}
For a path $p=(i_1,\dots,i_t)$ composed of a set of customers respecting the monotonicity condition, the expected recourse cost decreases monotonically when customers are removed from $p$, i.e., $\mathcal{Q}(p) \geq \mathcal{Q}(p')$ for any subpath $p'$ of $p$.
\end{proposition}
\proof{}
\begin{align*}
    \mathcal{Q}(p) &= 2\sum_{j=1}^t\sum_{l=1}^{\infty}\mathds{P}\left(\sum_{k=1}^{j-1}\xi_{i_k} \leq lQ < \sum_{k=1}^{j}\xi_{i_k}\right)c_{0i_j}\\
    &= 2\sum_{\substack{j=1 \\ i_j \in p'}}^t\sum_{l=1}^{\infty}\mathds{P}\left(\sum_{k=1}^{j-1}\xi_{i_k} \leq lQ < \sum_{k=1}^{j}\xi_{i_k}\right)c_{0i_j} + 2\sum_{\substack{j=1 \\ i_j \in p \setminus p'}}^t\sum_{l=1}^{\infty}\mathds{P}\left(\sum_{k=1}^{j-1}\xi_{i_k} \leq lQ < \sum_{k=1}^{j}\xi_{i_k}\right)c_{0i_j}\\
    &\geq 2\sum_{\substack{j=1 \\ i_j \in p'}}^t\sum_{l=1}^{\infty}\mathds{P}\left(\sum_{k=1}^{j-1}\xi_{i_k} \leq lQ < \sum_{k=1}^{j}\xi_{i_k}\right)c_{0i_j}\\
    &\geq 2\sum_{\substack{j=1 \\ i_j \in p'}}^t\sum_{l=1}^{\infty}\mathds{P}\left(\sum_{\substack{k=1 \\ i_k \in p'}}^{j-1}\xi_{i_k} \leq lQ < \sum_{\substack{k=1 \\ i_k \in p'}}^{j}\xi_{i_k}\right)c_{0i_j}\\
    &= \mathcal{Q}(p')
\end{align*}

The second inequality comes from the fact that, for any $l{\geq}1$ and for any $j \in \{1,\dots,t\}$ such that $i_j \in p'$, the following inequality holds.
\begin{equation}
    \label{monotonicity_condition_exp_rec} \mathds{P}\left(\sum_{k=1}^{j-1}\xi_{i_k} \leq lQ < \sum_{k=1}^{j}\xi_{i_k}\right) \geq \mathds{P}\left(\sum_{\substack{k=1 \\ i_k \in p'}}^{j-1}\xi_{i_k} \leq lQ < \sum_{\substack{k=1 \\ i_k \in p'}}^{j}\xi_{i_k}\right)
\end{equation}

To prove this, we first introduce some notation:
\begin{itemize}
    \item $\tilde{S}^j = \left\{i_k \in p' : k \leq j{-}1\right\}$ is the set of customers that are visited before customer $i_j$ in path $p'$
    \item $\tilde{\xi}^j = \sum_{i \in \tilde{S}^j}\xi_{i}$ is the total demand of these customers
    \item $\{a_1,\dots,a_s\} = \left\{i_k \in p : k \leq j{-}1, i_k \notin p' \right\}$ is the set of customers that are visited before customer $i_j$ in path $p$, but are not visited in path $p'$
\end{itemize}

We can now prove inequality (\ref{monotonicity_condition_exp_rec}).
\begingroup
\allowdisplaybreaks
\begin{align}
    \mathds{P}\left(\sum_{\substack{k=1 \\ i_k \in p'}}^{j-1}\xi_{i_k} \leq lQ < \sum_{\substack{k=1 \\ i_k \in p'}}^{j}\xi_{i_k}\right) &= \mathds{P}\left(\tilde{\xi}^j \leq lQ < \tilde{\xi}^j+\xi_{i_j}\right) \\
    \label{monotonic_expected_recourse_ineq1}&\leq \mathds{P}\left(\tilde{\xi}^j + \xi_{a_1} \leq lQ < \tilde{\xi}^j + \xi_{a_1} +\xi_{i_j}\right)\\
    \label{monotonic_expected_recourse_ineq2}& \ \vdots \\
    \label{monotonic_expected_recourse_ineq3}&\leq \mathds{P}\left(\tilde{\xi}^j + \xi_{a_1} + \dots + \xi_{a_s} \leq lQ < \tilde{\xi}^j  + \xi_{a_1} + \dots + \xi_{a_s} +\xi_{i_j} \right)\\
    & = \mathds{P}\left(\sum_{k=1}^{j-1}\xi_{i_k} \leq lQ < \sum_{k=1}^{j}\xi_{i_k}\right)
\end{align}
\endgroup

By hypothesis, the set $S=\{i_1, \dots, i_j\}$ respects the monotonicity condition. This allows us to apply inequality (\ref{ineq:monotonicity_condition_1}) $s$ times for different pairs $(a,b) \in S \times S$ and subsets $\tilde{S} \subseteq S$ to obtain inequalities (\ref{monotonic_expected_recourse_ineq1})-(\ref{monotonic_expected_recourse_ineq3}). At the $q$-th application of the inequality, $(a,b)=(a_q,i_j)$ and $\tilde{S}=\tilde{S}^j \cup \{a_1, \dots, a_{q-1}\}$.
\endproof

Propositions \ref{prop:poisson} and \ref{prop:normal} show that, under some conditions, the monotonicity property is respected for customer demands  following independent Poisson, normal, binomial, Erlang, and negative binomial distributions.

\begin{proposition}
\label{prop:poisson}
Any set $S$ of customers whose demands are given by independent Poisson RVs $\xi_i \sim \text{Poisson}(\lambda_i)$ with total expected value $\sum_{i\in S} \lambda_{i} \leq Q$ respects the monotonicity condition.
\end{proposition}

\proof{}
Consider a pair $(a,b) \in S \times S$ such that $a \neq b$, a subset $\tilde{S}\subseteq S \setminus \{a,b\}$ with total expected demand $\tilde{\lambda} = \sum_{i \in \tilde{S}}\lambda_i$ and a positive integer $l \in \mathds{N}$. To show that inequality (\ref{ineq:monotonicity_condition_2}) holds, we use the fact that the sum of $s$ independent Poisson RVs with means $\lambda_1, \dots, \lambda_s$ is Poisson distributed with mean $\sum_{k=1}^s \lambda_k$. Also, we use the fact that the cumulative distribution function (CDF) of a Poisson RV with mean $\lambda$ evaluated at $lQ \in \mathds{N}$ is given by $H(lQ,\lambda)=\frac{\Gamma(lQ{+}1, \lambda)}{(lQ)!}$, where $\Gamma(lQ{+}1, \lambda)=\int_{\lambda}^{\infty}t^{lQ}e^{-t}dt$ is the upper incomplete gamma function.
\begingroup
\allowdisplaybreaks
\begin{align}
    &\mathds{P}\left(\tilde{\xi} + \xi_a \leq lQ\right) - \mathds{P}\left(\tilde{\xi} + \xi_a + \xi_b \leq lQ\right) \geq \mathds{P}\left(\tilde{\xi} \leq lQ\right) - \mathds{P}\left(\tilde{\xi} + \xi_b \leq lQ\right)\\
    \iff & H\left(lQ, \tilde{\lambda}{+} \lambda_a\right) - H\left(lQ, \tilde{\lambda}{+}\lambda_a{+} \lambda_b\right) \geq H\left(lQ, \tilde{\lambda}\right) - H\left(lQ, \tilde{\lambda}{+}\lambda_b\right)\\
    \iff & \frac{\Gamma\left(lQ+1, \tilde{\lambda}{+} \lambda_a\right) }{(lQ)} - \frac{\Gamma\left(lQ+1, \tilde{\lambda}{+}\lambda_a{+} \lambda_b\right) }{(lQ)} \geq \frac{\Gamma\left(lQ+1, \tilde{\lambda}\right) }{(lQ)} - \frac{\Gamma\left(lQ+1, \tilde{\lambda}{+}\lambda_b\right) }{(lQ)} \\
    \iff & \Gamma\left(lQ+1, \tilde{\lambda}{+} \lambda_a\right)  - \Gamma\left(lQ+1, \tilde{\lambda}{+}\lambda_a{+} \lambda_b\right)  \geq \Gamma\left(lQ+1, \tilde{\lambda}\right)  - \Gamma\left( lQ+1 , \tilde{\lambda}{+}\lambda_b\right) \\
    \iff & \int_{ \tilde{\lambda}{+} \lambda_a }^{\infty} \lambda^{ lQ } e^{-\lambda} d\lambda - \int_{ \tilde{\lambda}{+}\lambda_a{+} \lambda_b }^{\infty} \lambda^{ lQ } e^{-\lambda} d\lambda \geq \int_{ \tilde{\lambda} }^{\infty} \lambda^{ lQ } e^{-\lambda} d\lambda - \int_{ \tilde{\lambda}{+}\lambda_b }^{\infty} \lambda^{ lQ } e^{-\lambda} d\lambda \\
    \iff & \int_{ \tilde{\lambda}{+} \lambda_a }^{\tilde{\lambda}{+}\lambda_a{+} \lambda_b} \lambda^{ lQ } e^{-\lambda} d\lambda  \geq \int_{ \tilde{\lambda} }^{\tilde{\lambda}{+}\lambda_b} \lambda^{ lQ } e^{-\lambda} d\lambda\\
    \iff & \int_{ \tilde{\lambda} }^{\tilde{\lambda}{+} \lambda_b} \left(\lambda{+}\lambda_a\right)^{ lQ } e^{-\left(\lambda{+}\lambda_a\right)} d\lambda  \geq \int_{ \tilde{\lambda} }^{\tilde{\lambda}{+}\lambda_b} \lambda^{ lQ } e^{-\lambda} d\lambda\\
    \iff & \int_{ \tilde{\lambda} }^{\tilde{\lambda}{+} \lambda_b} \left( \left(\lambda{+}\lambda_a\right)^{ lQ } e^{-\left(\lambda{+}\lambda_a\right)} - \lambda^{ lQ } e^{-\lambda} \right)d\lambda  \geq 0\\
    \label{poisson_last_ineq} \iff & \int_{ \tilde{\lambda} }^{\tilde{\lambda}{+} \lambda_b} \left( f_{ lQ }(\lambda + \lambda_a) - f_{ lQ }(\lambda) \right)d\lambda  \geq 0,
\end{align}
\endgroup
where $f_{ lQ }(\lambda) := \lambda^{ lQ }e^{-\lambda}$. Since the inequalities $0 \leq \lambda \leq \lambda{+}\lambda_{a} \leq  Q  \leq  lQ $ hold for any $\lambda \in [\tilde{\lambda}, \tilde{\lambda}{+}\lambda_b]$, it is sufficient to show that the function $f_{ lQ }(\cdot)$ is non-decreasing on $[0, lQ ]$ to conclude that the difference $f_{ lQ }(\lambda + \lambda_a) - f_{ lQ }(\lambda)$ is non-negative everywhere on the interval of integration. The derivative of function $f_{ lQ }(\cdot)$ is given by:
\begin{align*}
    \frac{\partial f_{ lQ }(\lambda)}{\partial \lambda} &= ({ lQ }-\lambda)e^{-\lambda}\lambda^{{ lQ }-1},
\end{align*}
which is indeed non-negative $\forall$ $\lambda \in [0, lQ ]$. Inequality (\ref{poisson_last_ineq}) is therefore respected.
\endproof

\begin{proposition}
\label{prop:normal}
Any set $S$ of customers whose demands are given by independent normal RVs $\xi_i \sim \mathcal{N}(\mu_i,\sigma^2_{i})$ with integer expected values $\mu_i \in \mathds{N}$ and sharing a common coefficient of dispersion $D=\frac{\sigma^2_{i}}{\mu_i}\leq 1$ $\forall i \in S$ with total expected value $\sum_{i\in S} \mu_{i} \leq Q$ respects the monotonicity condition.
\end{proposition}
\proof{}
Although the support of a normally distributed RV is $\mathds{R}$, the demands are regarded as non-negative in the proof. For this approximation not to be problematic in practice, the standard deviation $\sigma_i$ of each customer $i \in N$ should not exceed a small fraction of its demand $\mu_i$.

Consider a pair $(a,b) \in S \times S$ such that $a \neq b$ and a subset $\tilde{S}\subseteq S \setminus \{a,b\}$ with total expected demand $\tilde{\mu}=\sum_{i \in \tilde{S}}\mu_i$. Using the notation of Proposition \ref{prop:monotonicity_condition_sufficient}, $\tilde{\xi}$, $\xi_a$ and $\xi_b$ can respectively be decomposed as the sum of $\tilde{\rho}{=}\tilde{\mu}$, $\rho_a{=}\mu_a$ and $\rho_b{=}\mu_b$ i.i.d. RVs $\xi_k$, $k \in \mathds{N}$, such that $\zeta_k \sim \mathcal{N}(\mu{=}1, \sigma^2{=}D)$. From there, the sufficient condition (\ref{ineq:monotonicity_condition_sufficient}) can be rewritten as follows.
\begin{equation}
\label{ineq:monotonicity_condition_sufficient_normal} P_{lQ}(\tilde{\mu}+\mu_a+j) \geq P_{lQ}(\tilde{\mu}+j) \hspace{0.5cm} \forall l \in \mathds{N}, \forall j \in \{1,\dots,\mu_b\}
\end{equation}

The coefficient of variation $c=\sigma/\mu$ of the i.i.d. RVs $\zeta_k$ is equal to $\sqrt{D}$. Since $D \leq 1$ by hypothesis, then $c\leq 1$. In this case, it follows from \cite{KREIMER199063} that, for any $l \in \mathds{N}$, the sequence $\{P_{lQ}(s)\}_{s=1}^{\lfloor lQ \rfloor}$ is monotonically increasing. Since $\tilde{\mu}{+}\mu_a{+}\mu_b \leq \sum_{i \in S}\mu_i \leq Q$ and all the expected demands are integer-valued, this allows to conclude that (\ref{ineq:monotonicity_condition_sufficient_normal}) holds.
\endproof

\begin{proposition}
\label{prop:binomial}
Any set $S$ of customers whose demands are given by independent binomial RVs $\xi_i \sim \text{Bin}(n_i,p)$ sharing a common success probability $p$ with total expected value $\sum_{i\in S} n_{i}p \leq Q$ respects the monotonicity condition.
\end{proposition}
\proof{}
Consider a pair $(a,b) \in S \times S$ such that $a \neq b$ and a subset $\tilde{S}\subseteq S \setminus \{a,b\}$ with total expected demand $\tilde{n}p$, where $\tilde{n}=\sum_{i \in \tilde{S}}n_i$. Using the notation of Proposition \ref{prop:monotonicity_condition_sufficient}, the demands $\tilde{\xi}$, $\xi_a$ and $\xi_b$ can respectively be decomposed as the sum of $\tilde{\rho}{=}\tilde{n}$, $\rho_a{=}n_a$ and $\rho_b{=}n_b$ i.i.d. RVs $\xi_k$, $k \in \mathds{N}$, such that $\zeta_k \sim \text{Bern}(p)$. From there, the sufficient condition (\ref{ineq:monotonicity_condition_sufficient}) can be rewritten as follows.
\begin{equation}
\label{ineq:monotonicity_condition_sufficient_binomial} P_{lQ}(\tilde{n}+n_a+j) \geq P_{lQ}(\tilde{n}+j) \hspace{0.5cm} \forall l \in \mathds{N}, \forall j \in \{1,\dots,n_b\}
\end{equation}

It follows from \cite{KREIMER199063} that, for any $l \in \mathds{N}$, the sequence $\{P_{lQ}(s)\}_{s=1}^{\left\lceil \frac{lQ}{p} \right\rceil }$ is non-decreasing. Since $\tilde{n}{+}n_a{+}n_b \leq \sum_{i\in S} n_{i} \leq Q/p$, this allows to conclude that (\ref{ineq:monotonicity_condition_sufficient_binomial}) holds.
\endproof

\begin{proposition}
\label{prop:erlang}
Any set $S$ of customers whose demands are given by independent Erlang RVs $\xi_i \sim \text{Erlang}(n_i,\lambda)$ sharing a common rate parameter $\lambda$ with total expected value $\sum_{i\in S} n_i/\lambda \leq Q$ respects the monotonicity condition.
\end{proposition}
\proof{}
Consider a pair $(a,b) \in S \times S$ such that $a \neq b$ and a subset $\tilde{S}\subseteq S \setminus \{a,b\}$ with total expected demand $\tilde{n}/\lambda$, where $\tilde{n}=\sum_{i \in \tilde{S}}n_i$. Using the notation of Proposition \ref{prop:monotonicity_condition_sufficient}, the demands $\tilde{\xi}$, $\xi_a$ and $\xi_b$ can respectively be decomposed as the sum of $\tilde{\rho}{=}\tilde{n}$, $\rho_a{=}n_a$ and $\rho_b{=}n_b$ i.i.d. RVs $\xi_k$, $k \in \mathds{N}$, such that $\zeta_k \sim \text{Exp}(\lambda)$. From there, the sufficient condition (\ref{ineq:monotonicity_condition_sufficient}) can be rewritten as follows.
\begin{equation}
\label{ineq:monotonicity_condition_sufficient_erlang} P_{lQ}(\tilde{n}+n_a+j) \geq P_{lQ}(\tilde{n}+j) \hspace{0.5cm} \forall l \in \mathds{N}, \forall j \in \{1,\dots,n_b\}
\end{equation}

It follows from \cite{KREIMER199063} that, for any $l \in \mathds{N}$, the sequence $\{P_{lQ}(s)\}_{s=1}^{\left\lceil \lambda lQ \right\rceil }$ is non-decreasing. Since $\tilde{n}{+}n_a{+}n_b \leq \sum_{i\in S} n_{i} \leq \lambda Q$, this allows to conclude that (\ref{ineq:monotonicity_condition_sufficient_erlang}) holds.
\endproof

\begin{proposition}
\label{prop:negative_binomial}
Any set $S$ of customers whose demands are given by independent negative binomial RVs $\xi_i \sim \text{NB}(r_i,p)$ sharing a common success probability $p$ with total expected value $\sum_{i\in S} r_{i}(1{-}p)/p \leq Q$ respects the monotonicity condition.
\end{proposition}
\proof{}
Consider a pair $(a,b) \in S \times S$ such that $a \neq b$, a subset $\tilde{S}\subseteq S \setminus \{a,b\}$ with total expected demand $\tilde{r}(1{-}p)/p$, where $\tilde{r}=\sum_{i \in \tilde{S}}r_i$ and a positive integer $l \in \mathds{N}$. We use the fact that the sum of $s$ independent negative binomial RVs with parameters $r_1, \dots, r_s$ and the same success probability $p$ is also a negative binomial RV with parameters $r=\sum_{k=1}^s r_k$ and $p$ to rewrite inequality (\ref{ineq:monotonicity_condition_1}) as follows.
\begin{align}
    &\mathds{P}\left(\tilde{\xi} + \xi_a \leq lQ <\tilde{\xi} + \xi_a + \xi_b \right) \geq \mathds{P}\left(\tilde{\xi} \leq lQ <\tilde{\xi} + \xi_b\right)\\
    \label{neg_bin_eq1}\iff& \mathds{P}\left(\sum_{k=1}^{lQ+\tilde{r}+r_a} \zeta_k \leq lQ < \sum_{k=1}^{lQ+\tilde{r}+r_a+r_b} \zeta_k  \right) \geq \mathds{P}\left(\sum_{k=1}^{lQ+\tilde{r}} \zeta_k \leq lQ < \sum_{k=1}^{lQ+\tilde{r}+r_b} \zeta_k  \right)
\end{align}

To obtain inequality (\ref{neg_bin_eq1}), we use the fact that the CDF of a negative binomial with parameters $r$ and $p$ evaluated at $lQ \in \mathds{N}$ is equal to the CDF of a binomial distribution with $lQ+r$ trials and probability of success $1-p$ evaluated at $lQ$. In inequality (\ref{neg_bin_eq1}), these binomial RVs are decomposed into sums of i.i.d. RVs $\zeta_k \sim \text{Bern}(1{-}p)$. Using this reformulation of the monotonicity condition based on the RVs $\zeta_k$, it suffices to show that $P_{lQ}(lQ+\tilde{r}+r_a+j) \geq P_{lQ}(lQ+\tilde{r}+j) \ \forall j \in \{1,\dots,r_b\}$ to conclude the proof, by Proposition \ref{prop:monotonicity_condition_sufficient}. For a given $j\in \{1,\dots,r_b\}$, this inequality can be developed as follows:
\begingroup
\allowdisplaybreaks
\begin{align}
    &P_{lQ}(lQ+\tilde{r}+r_a+j) \geq P_{lQ}(lQ+\tilde{r}+j)\\
    \iff& \mathds{P}\left( \sum_{k=1}^{lQ+\tilde{r}+r_a+j-1} \zeta_k \leq lQ < \sum_{k=1}^{lQ+\tilde{r}+r_a+j} \zeta_k \right) \geq \mathds{P}\left( \sum_{k=1}^{lQ+\tilde{r}+j-1} \zeta_k \leq lQ < \sum_{k=1}^{lQ+\tilde{r}+j} \zeta_k \right)\\
    \label{neg_bin_eq2}\iff& \mathds{P}\left( \chi \leq lQ < \chi + \zeta_1 \right) \geq \mathds{P}\left(\bar{\chi} \leq lQ < \bar{\chi} + \zeta_1 \right),
\end{align}
\endgroup
where $\chi\sim \text{Bin}(lQ+\tilde{r}{+}r_a{+}j{-}1, 1{-}p)$ and $\bar{\chi} \sim \text{Bin}(lQ{+}\tilde{r}{+}j{-}1, 1{-}p)$. Note that the event in the left-hand side probability of inequality (\ref{neg_bin_eq2}) occurs if and only if $\chi{=}lQ$ and $\zeta_1{=}1$. Since $\chi$ and $\zeta_1$ are independent, the joint probability of these two events is equal to the product of their probabilities. Applying the same reasoning to the left-hand side, inequality (\ref{neg_bin_eq2}) can be written as follows.
\begingroup
\allowdisplaybreaks
\begin{align}
    &\mathds{P}\left( \chi = lQ \right)\mathds{P}\left( \zeta_1 = 1 \right) \geq \mathds{P}\left( \bar{\chi} = lQ \right)\mathds{P}\left( \zeta_1 = 1 \right)\\
    \iff& \mathds{P}\left( \chi = lQ \right) \geq \mathds{P}\left( \bar{\chi} = lQ \right)\\
    \iff& \frac{(lQ+\tilde{r}+r_a+j-1)!}{(lQ)!(\tilde{r}+r_a+j-1)!}(1-p)^{lQ}p^{\tilde{r}+r_a+j-1} \geq \frac{(lQ+\tilde{r}+j-1)!}{(lQ)!(\tilde{r}+j-1)!}(1-p)^{lQ}p^{\tilde{r}+j-1}\\
    \iff& \frac{\prod_{i=1}^{r_a}(lQ+\tilde{r}+i+j-1)}{\prod_{i=1}^{r_a}(\tilde{r}+i+j-1)}p^{r_a} \geq 1\\
    \label{neg_bin_eq3}\iff& \prod_{i=1}^{r_a}\left(\frac{p(lQ+\tilde{r}+i+j-1)}{\tilde{r}+i+j-1}\right) \geq 1
\end{align}
\endgroup

To complete the proof, we demonstrate that inequality (\ref{neg_bin_eq3}) holds by showing that each term in the product is greater than or equal to 1. For $i \in \{1,\dots,r_a\}$, we have:
\begingroup
\allowdisplaybreaks
\begin{align}
    \frac{p(lQ+\tilde{r}+i+j-1)}{\tilde{r}+i+j-1} &= \frac{plQ}{\tilde{r}+i+j-1}+p\\
    \label{neg_bin_eq4}&\geq \frac{pQ}{\tilde{r}+r_a+r_b}+p\\
    &\geq \frac{pQ}{\sum_{i \in S}r_i}+p\\
    \label{neg_bin_eq5}&\geq \frac{(1-p)pQ}{pQ}+p\\
    &= 1
\end{align}
\endgroup
To obtain inequality (\ref{neg_bin_eq4}), we notice that $l\geq 1$ and $i{+}j{-}1 \leq r_a{+}r_b{-}1 \leq r_a{+}r_b$. Inequality (\ref{neg_bin_eq5}) follows from the assumption that the total expected demand of the customers of set $S$ does not exceed the vehicle capacity.
\endproof

The next propositions present two counterexamples to illustrate that the recourse function is not monotonic in general.

\begin{proposition}
\label{prop:capacity constraints}
Let $p=(i_1,\dots,i_t)$ be a path of customers whose demands are given by independent Poisson RVs $\xi_i \sim \text{Poisson}(\lambda_i)$ with total expected value $\sum_{i=1}^{t}\lambda_i > Q$, and $p'$ a subpath of $p$. Inequality $\mathcal{Q}(p) \geq \mathcal{Q}(p')$ does not hold in general.
\end{proposition}
\proof{}
Consider a problem with $Q=20$, a path $p=(1,2,3)$, and $p'=(2,3)$. Suppose the expected demand of each customer is given by $\lambda_1=5$, $\lambda_2=15$, and $\lambda_3=10$, respectively. If $c_{01}=c_{02}=0$ and $c_{03}=1$, then $\mathcal{Q}(p) \approx 1.11 < 1.47 \approx \mathcal{Q}(p')$.
\endproof

\begin{proposition}\label{prop:Non monotonic}
Let $p=(i_1,\dots,i_t)$ be a path of customers whose total expected demand respects the vehicle capacity and $p'$ a subpath of $p$. Inequality $\mathcal{Q}(p) \geq \mathcal{Q}(p')$ does not hold in general.
\end{proposition}
\proof{}
Consider a problem with $Q=20$, a path $p=(1,2,3)$, and $p'=(2,3)$. Suppose the demand of customer 1 is 5 with probability 1, while that of customers 2 and 3 is 6 with probability 0.9 and 16 with probability 0.1. If $c_{03} > c_{02}$, then $\mathcal{Q}(p) = (0.1)2c_{02} + (0.09)2c_{03} < (0.19)2c_{03} = \mathcal{Q}(p')$.
\endproof

\section{Disaggregated integer L-shaped method} \label{sec:Disaggregated integer L-shaped method}

This section presents the DL-shaped method as a variant of the classical integer L-shaped method from \cite{laporte1993integer}. First, an outline of the classical method is provided, and we motivate the need for our method by exposing the known shortcoming of the classical one. Next, Section \ref{sec:New Opt Cuts} presents new optimality cuts, from which the validity of our method results. Finally, Section \ref{sec:LBFs} introduces the new LBFs we use to approximate the recourse cost.

The integer L-shaped method for two-stage recourse problems starts by constructing a so-called master problem (MP). The MP is obtained from the original problem by relaxing a set of complicating constraints. This set typically includes the integrality constraints and any set of constraints that is exponential in size, like the classical subtour elimination constraints. In addition, the recourse function is replaced by a continuous variable $\theta \geq 0$ that bounds it from below. If a lower bound on the expected recourse $L$ can be computed, then the inequality $\theta \geq L$ is added to the model. Otherwise, the bound on $\theta$ remains the trivial one. The MP is thus a problem defined by the non-relaxed constraints of the original program, the non-negativity constraint on $\theta$, the lower bound $L$ on $\theta$, and an objective function of the form $\min c x + \theta$, where $c$ is a known cost vector with the same size as the solution vector. The problem is then solved with an iterative branch-and-bound method that adds linear inequalities on the fly to explore the set of feasible solutions of the original program. The first MP that is solved in the method is denoted as the root node, and the following MPs are indexed by the iteration at which they appear in the resolution process. The linear inequalities that are added to the successive MPs are either feasibility or optimality cuts. The feasibility cuts prohibit infeasible solutions from being visited, while optimality cuts improve the bound on $\theta$ whenever a better feasible solution is found. The feasibility cuts are generally problem-specific. In the context of the \vrpsd, they typically correspond to subtour elimination constraints and rounded capacity inequalities. Optimality cuts have a general form regardless of the problem to which the integer L-shaped method is applied. Let $x^\nu_1$ be the set of positive variables in a feasible solution $x^\nu$ at iteration $\nu$ of the method. Then the following inequality:
\begin{equation}
  \label{eq:optcut1993} \theta  \geq (\mathcal{Q}(x^\nu) - L)\left ( \sum_{i \in x^\nu_1} x_i - \sum_{i \not\in x^\nu_1} x_i - |x^\nu_1| + 1 \right ) + L,
\end{equation}
where $\mathcal{Q}(x^\nu)$ is the expected recourse cost of solution $x^\nu$, is a valid optimality cut. It was shown by \cite{laporte1993integer} that these optimality cuts will yield an optimal solution, if one exists, in a finite number of steps.

The known shortcoming of the integer L-shaped method is that its optimality cuts are active for a single solution. The optimality cut reduces to $\theta \geq \mathcal{Q}(x^\nu)$ at solution $x^\nu$, but for any other solution, it does not improve the inequality $\theta \geq L$. Consequently, without appropriate LBFs to increase the bound on $\theta$, the classical method generally iterates through many first-stage solutions before converging to an optimal solution.

\subsection{New optimality cuts}\label{sec:New Opt Cuts}

The main idea of the DL-shaped method is to replace the variable $\theta$ that bounds the recourse function by a sum of variables $\theta_i$, one for each customer $i \in N$. The purpose of this transformation is to express the contribution of each customer to the cost of the second stage by the gradual addition of a new type of optimality cuts during the resolution. These cuts are added for feasible paths that are found, and the cut of each path only involves the variables associated with the arcs of the path and the $\theta_i$ variables of the customer of the path. This makes them active for each path that is active in the current solution, whether it is fractional or integer. They allow to bound the recourse function more tightly and effectively than the traditional optimality cuts, as those only bound one solution per cut. This gives the following objective function.


\begin{equation}
\min \sum_{(i,j)\in E}c_{ij}x_{ij} + \sum_{i \in N} \theta_i \label{model1:vrp7}
\end{equation}
The DL-shaped method starts by constructing a MP with the new objective function, relaxed integrality constraints, and rounded capacity inequalities. In each iteration of the method, when a feasible integer solution is found, an optimality cut is added for each route $r=(0,p,0)$ having a positive recourse ($\mathcal{Q}(r)>0$), where $p=(i_1,\ldots,i_t)$. Such a path $p$ is composed of $|p|=t-1$ edges. For a first-stage solution $x$, let us denote by $x(p)=\sum_{j=1}^{t-1}x_{i_{j}i_{j+1}}$ the number of these edges that are active. Also, let $N(r)$ be the set of customers that are visited in path $p$. The new optimality cut follows.
\begin{equation}
\sum_{i \in N(r)} \theta_i \geq \mathcal{Q}(r)\left( x(p) - |p|+1 \right)  \label{optcut_path}
\end{equation}

We refer to these cuts as path cuts (P-Cuts). Since they do not include the depot edges, as opposed to the route cuts of \cite{seguin1996problemes}, these P-Cuts are more general. In the specific case of a route $r=(0,i,0)$ visiting a single customer $i$, inequality (\ref{optcut_path}) reduces to $\theta_i \geq \mathcal{Q}(r)$.

We now prove that the P-Cuts are valid optimality cuts. By \cite{laporte1993integer}, this implies that the DL-shaped method yields an optimal solution in a finite number of steps.

\begin{proposition}
Let $x^\nu$ be a first-stage feasible solution and $\mathcal{Q}(x^\nu)$ the corresponding recourse cost. The set of optimality cuts (\ref{optcut_path}), for $r \in \mathcal{R}^{\nu}$, is valid if the problem has the monotonicity property.\label{prop:PCuts validity}
\end{proposition}
\proof{}
Let $x^{\lambda}$ be a solution satisfying constraints (\ref{model:depot_degree})-(\ref{model:int}). For the right-hand side of (\ref{optcut_path}) we have $\left( x^{\lambda}(p) - |p|+1 \right) \leq 1$, with equality if and only if path $p$ is a connected subpath of a route $r$ contained in solution $x^{\lambda}$. In particular, for $x^{\lambda} = x^\nu$, for each route $r \in \mathcal{R}^{\nu}$, the right-hand side of (\ref{optcut_path}) takes the value $\mathcal{Q}(r)$. Together, these $|\mathcal{R}^{\nu}|$ constraints imply that a solution $(x^\nu,\theta^\nu)$ respects $\sum_{i \in N} \theta^{\nu}_i = \sum_{r \in \mathcal{R}^{\nu}}\sum_{i \in N(r)} \theta^{\nu}_i \geq \sum_{r \in \mathcal{R}^{\nu}} \mathcal{Q}(r) = \mathcal{Q}(x^\nu)$. Therefore, these cuts bind $\sum_{i \in N} \theta^{\nu}_i$ to $\mathcal{Q}(x^\nu)$ and the objective value of a feasible solution $(x^\nu,\theta^\nu)$ for the MP cannot be smaller than that of $x^\nu$ for problem (\ref{model:obj})-(\ref{model:int}).

Let $\theta^*(x^{\lambda})$ be the minimizer of $\sum_{i \in N}\theta^{\lambda}_i$ over the feasible set specified by the first-stage solution $x^{\lambda}$ and by the constraints (\ref{optcut_path}) that are contained in the MP at a given iteration. Next, we prove that constraints (\ref{optcut_path}) produce valid lower bounds on $\sum_{i \in N}\theta^{\lambda}_i$. To do so, we show that the objective value of a solution $(x^{\lambda},\theta^*(x^{\lambda}))$ to the MP is never larger than that of $x^{\lambda}$ for problem (\ref{model:obj})-(\ref{model:int}). The MP is maximally constrained if (\ref{optcut_path}) has been added to the model for each possible route $r$. Let $\mathcal{R^{\lambda}}$ be the set of routes of solution $x^{\lambda}$. Under this maximal set of optimality cuts, for each $r \in \mathcal{R^{\lambda}}$, the sum $\sum_{i \in N(r)}\theta^{\lambda}_i$ is constrained to be larger than or equal to $F(\tilde{r})$ for each route $\tilde{r}=(0,\tilde{p},0)$ such that $\tilde{p}$ is a connected subpath of route $r$. Since it is assumed that the problem respects the monotonicity property, it follows from Proposition \ref{prop:monotonic_expected_recourse} that the most restrictive of these constraints is $\sum_{i \in N(r)}\theta^{\lambda}_i \geq \mathcal{Q}(r)$. Hence, by definition of $\theta^*(x^{\lambda})$, $\sum_{i \in N(r)}\theta_i^*(x^{\lambda})=\mathcal{Q}(r)$ and thus, $\sum_{i \in N} \theta^*_i(x^{\lambda}) = \sum_{r \in \mathcal{R}^{\lambda}}\sum_{i \in N(r)} \theta^*_i(x^{\lambda}) = \sum_{r \in \mathcal{R}^{\lambda}} \mathcal{Q}(r) = \mathcal{Q}(x^{\lambda})$. The objective value of $(x^{\lambda},\theta^*(x^{\lambda}))$ and $x^{\lambda}$ for their respective problems are therefore identical when the MP is maximally constrained. When cuts are still to be generated, the objective value of $(x^{\lambda},\theta^*(x^{\lambda}))$ for the MP can thus only be less than or equal to that of $x^{\lambda}$ for problem (\ref{model:obj})-(\ref{model:int}).
\endproof

As shown in \cite{laporte1993integer}, the linear relaxation of \eqref{model:obj}-\eqref{model:int} can be improved by introducing a general lower bound on the recourse using a precomputed value $L$.
\begin{equation}
    \sum_{i \in N}\theta_i \geq L \label{eq:L1}
\end{equation}

Unfortunately, such a bound is generally close to 0 for problems with an unlimited fleet, as using more vehicles reduces the probability of having a failure. In particular, the recourse cost is close to 0 when each vehicle only serves a single client. It can therefore be useful to compute a lower bound $L_m$ on the recourse given that exactly $m \in M$ vehicles are being used. $L_m$ is computed using the method defined in Section \ref{sec:spf}. This leads to the following valid lower bound inequalities.
\begin{equation}
    \sum_{i \in N}\theta_i \geq \sum_{m\in M}L_m z_m \label{eq:L}
\end{equation}

To conclude this section, we propose a different optimality cut that does not require the recourse to be monotonic. In the absence of monotonicity, the validity of our method can be preserved by adding the depot edges to cut \eqref{optcut_path}. This way, the cut associated with a route $r=(0,i_1,\dots,i_t,0)$ is only active when this route appears in a first-stage solution. The coefficients of the non-depot edge variables are equal to two in order to avoid the cut being active in other solutions. The modified optimality cut follows.
\begin{equation}
    \sum_{i \in N(r)} \theta_i \geq \mathcal{Q}(r)\left( \sum_{j=1}^{t-1} 2x_{i_ji_{j+1}} + x_{0i_1} + x_{0i_t} - 2(t-1)-1 \right) \label{eq:RCuts}
\end{equation}

\subsection{Lower bounding functionals}\label{sec:LBFs}

This section presents a new type of LBFs that extends some ideas proposed in \cite{cote2020vehicle}. Similarly to the previously introduced optimality cuts, this type of cuts uses variables to bound the recourse for sets of customers $S$. The LBFs of \cite{cote2020vehicle} use a different set of variables than those used for optimality cuts, which has the disadvantage of making the LBFs independent from the optimality cuts. This work improves over them by using the same $\theta_i$ variables for both types of cuts. The new LBFs bound the recourse using all the $\theta_i$ variables associated with a set of customers $S$. A cut is active whenever a path visits all the customers of this set consecutively. The new LBFs are referred to as the set cuts (S-Cuts) and are as follows:
\begin{equation}
\sum_{i \in S}\theta_i \geq L(S)\left(x(S) - |S| +  \left\lceil \frac{\sum_{i \in S}\mu_i}{Q} \right\rceil +1 \right)  \label{eq:Set Cuts},
\end{equation}
where $x(S)=\sum_{(i,j) \in E(S)}x_{ij}$. The   S-Cut (\ref{eq:Set Cuts}) of a set $S$ is based on the least possible recourse cost of visiting the customers of $S$. Such a lower bound $L(S)$ is valid if, for any set of feasible paths $\{p_1, \dots, p_{m}\}$ forming a partition of $S$, the following inequality is respected.
\begin{equation}
    \label{ineq:L(S)} L(S) \leq \sum_{j = 1}^{m} \mathcal{Q}(p_j)
\end{equation}

\begin{proposition}
The S-Cut (\ref{eq:Set Cuts}) is valid if the problem has the monotonicity property.
\end{proposition}
\proof{}
Let $x^{\lambda}$ be a solution satisfying constraints (\ref{model:depot_degree})-(\ref{model:int}). If $x^{\lambda}(S) \leq |S|-\left\lceil \frac{\sum_{i \in S}\mu_i}{Q} \right\rceil - 1$, then the right-hand side of (\ref{eq:Set Cuts}) is non-positive and the constraint is thus respected by any feasible non-negative solution $(x^{\lambda}, \theta^{\lambda})$ to the MP. Otherwise, $x^{\lambda}(S) =  |S|-\left\lceil \frac{\sum_{i \in S}\mu_i}{Q} \right\rceil$. This means that solution $x^{\lambda}$ uses the smallest number of vehicles allowed by the capacity constraints to visit the customers of $S$ and that, in each route $r \in \mathcal{R}^{\lambda}$, the customers of $S$ are visited consecutively. In this case, the S-Cut reduces to $\sum_{i \in S}\theta_i\geq L(S)$. Let $\{p^{\lambda}_1, \dots, p^{\lambda}_{m}\}$ be the set of paths forming solution $x^{\lambda}$, where $m=\left\lceil \frac{\sum_{i \in S}\mu_i}{Q} \right\rceil$. To prove the validity of constraint (\ref{eq:Set Cuts}), we have to demonstrate that $L(S) \leq \sum_{j = 1}^{m} \mathcal{Q}(p^{\lambda}_j)$. For each $j \in \{1,\dots,m\}$, the segment of path $p^{\lambda}_j$ that is composed of customers of $S$ is a connected subpath $p_j$ of $p^{\lambda}_j$ and the set $\{p_1, \dots, p_{m}\}$ of these subpaths form a partition of $S$. Since it is assumed that the problem respects the monotonicity property, Proposition \ref{prop:monotonic_expected_recourse} implies that $\mathcal{Q}(p^{\lambda}_j) \geq \mathcal{Q}(p_j)$ $\forall j \in \{1,\dots,m\}$. By definition of the lower bound $L(S)$, the sequence of inequalities $ L(S) \leq \sum_{j = 1}^{m} \mathcal{Q}(p_j) \leq \sum_{j = 1}^{m} \mathcal{Q}(p^{\lambda}_j)$ allows to conclude.
\endproof

\subsection{Recourse lower bounds}\label{sec:Lower bounds}

This section presents lower bounds on the recourse used by the $S$-Cuts for a given set of customers $S$ to be served by $m$ vehicles. Some lower bounds from the literature are first presented.

To our knowledge, the first lower bound on the recourse of a given set $S$ is due to \cite{laporte2002integer}, where they suppose vehicles are only allowed to fail once. In the best case, the failures would happen for the customers of $S$ that are the closest to the depot. Then, the sum of demands is seen as a continuous quantity that has to be distributed among the $m$ vehicles. The problem is modeled as a non-linear program having for objective the minimization of the expected recourse.

A simpler bound, which is based on the assumption that there is one unique large vehicle with capacity $mQ$, is proposed by \cite{louveaux2018exact}. In the best case, the $l$-th failure would happen at the $l$-th customer closest to the depot. Then, the expected recourse cost can be easily calculated by summing the probability of the $l$-th failure multiplied by twice the $l$-th closest distance to the depot among the customers of $S$.

In the following, we propose three new bounds, $L_1(S)$, $L_2(S)$ and $L_3(S)$, which have the advantage of providing tighter values than those of \cite{laporte2002integer} and \cite{louveaux2018exact}. All these bounds are compared later in Section \ref{sec:Numerical results}.

\subsubsection{Single-route lower bound.} \label{sec:L for single-route}

The lower bound $L_1(S)$ introduced in this section applies to sets of customers that respect the monotonicity condition and can be assigned to a single vehicle. Let $\mathcal{P}^1(S)$ denote the set of all possible paths that visit the customers of $S$ consecutively. Proposition \ref{prop:lower_bound} shows that sorting these customers in non-increasing order of their distance to the depot produces the path of $\mathcal{P}^1(S)$ that minimizes the recourse cost. By the monotonicity condition, the probability that a restocking trip occurs at a given customer increases with its position in the path. The idea of this bound is thus to place customers that are far from the depot at the beginning of the path to avoid costly restocking trips.

\begin{proposition}
\label{prop:lower_bound}
Let $S \subseteq N$ be a set of customers respecting the monotonicity condition. Let $p^*=(i^*_1, i^*_2, \dots, i^*_{|S|}) \in \mathcal{P}^1(S)$ be a path such that $c_{0i^*_1} \geq c_{0i^*_2} \geq \dots \geq c_{0i^*_{|S|}}$. The expected recourse $\mathcal{Q}(p^*)$ of path $p^*$ is a valid lower bound on the expected recourse $\mathcal{Q}(p)$ of any path $p \in \mathcal{P}^1(S)$.
\end{proposition}
\proof{}
Let us consider a path $p = (i_1, \dots, i_{|S|}) \in \mathcal{P}^1(S)$ in which the customers are not sorted in non-increasing order of distances to the depot. In particular, this implies that there exists at least one index $h \in \left\{1,\dots, |S|{-}1\right\}$ such that $c_{0i_h} < c_{0i_{h+1}}$. We show that the expected recourse of path $p' = (i_1, \dots, i_{h{+}1}, i_h, \dots, i_{|S|}) \in \mathcal{P}^1(S)$, in which the order of visit of customers $i_{h}$ and $i_{h{+}1}$ is inverted, is less than or equal to that of path $p$.

We first develop the expression of $\mathcal{Q}(p)$:
\begingroup
\allowdisplaybreaks
\begin{align*}
    \mathcal{Q}(p) &= 2\sum_{j=1}^{|S|}\sum_{l=1}^{\infty}\mathds{P}\left(\sum_{k=1}^{j-1}\xi_{i_k} \leq lQ < \sum_{k=1}^{j}\xi_{i_k}\right)c_{0i_j} \\
    &= 2\sum\limits_{j=1}^{h-1}\sum\limits_{l=1}^{\infty}\mathds{P}\left(\sum\limits_{k=1}^{j-1}\xi_{i_k} \leq lQ < \sum\limits_{k=1}^{j}\xi_{i_k}\right)c_{0i_j} + 2\sum\limits_{l=1}^{\infty}\mathds{P}\left(\sum\limits_{k=1}^{h-1}\xi_{i_k} \leq lQ < \sum\limits_{k=1}^{h}\xi_{i_k}\right)c_{0i_h}\\
    &\ \ \ + 2\sum\limits_{l=1}^{\infty}\mathds{P}\left(\sum\limits_{k=1}^{h}\xi_{i_k} \leq lQ < \sum\limits_{k=1}^{h+1}\xi_{i_k}\right)c_{0i_{h+1}} + 2\sum\limits_{j=h+2}^{|S|}\sum\limits_{l=1}^{\infty}\mathds{P}\left(\sum\limits_{k=1}^{j-1}\xi_{i_k} \leq lQ < \sum\limits_{k=1}^{j}\xi_{i_k}\right)c_{0i_j}\\
    &= 2\sum\limits_{j=1}^{h-1}\sum\limits_{l=1}^{\infty}\mathds{P}\left(\sum\limits_{k=1}^{j-1}\xi_{i_k} \leq lQ < \sum\limits_{k=1}^{j}\xi_{i_k}\right)c_{0i_j} + 2\sum\limits_{l=1}^{\infty}\mathds{P}\left(\tilde{\xi} \leq lQ < \tilde{\xi} + \xi_{i_h}\right)c_{0i_h}\\
    &\ \ \ + 2\sum\limits_{l=1}^{\infty}\mathds{P}\left(\tilde{\xi} + \xi_{i_h} \leq lQ < \tilde{\xi} + \xi_{i_h} + \xi_{i_{h+1}} \right)c_{0i_{h+1}} + 2\sum\limits_{j=h+2}^{|S|}\sum\limits_{l=1}^{\infty}\mathds{P}\left(\sum\limits_{k=1}^{j-1}\xi_{i_k} \leq lQ < \sum\limits_{k=1}^{j}\xi_{i_k}\right)c_{0i_j},
\end{align*}
\endgroup
where $\tilde{\xi} = \sum_{k=1}^{h-1}\xi_{i_k}$. Doing the same for $\mathcal{Q}(p')$, we obtain:
\begingroup
\allowdisplaybreaks
\begin{align*}
    \mathcal{Q}(p') &= 2\sum\limits_{j=1}^{h-1}\sum\limits_{l=1}^{\infty}\mathds{P}\left(\sum\limits_{k=1}^{j-1}\xi_{i_k} \leq lQ < \sum\limits_{k=1}^{j}\xi_{i_k}\right)c_{0i_j} + 2\sum\limits_{l=1}^{\infty}\mathds{P}\left(\tilde{\xi} \leq lQ < \tilde{\xi} + \xi_{i_{h+1}}\right)c_{0i_{h+1}}\\
    &\ \ \ + 2\sum\limits_{l=1}^{\infty}\mathds{P}\left(\tilde{\xi} + \xi_{i_{h+1}} \leq lQ < \tilde{\xi} + \xi_{i_{h+1}} + \xi_{i_h} \right)c_{0i_h} + 2\sum\limits_{j=h+2}^{|S|}\sum\limits_{l=1}^{\infty}\mathds{P}\left(\sum\limits_{k=1}^{j-1}\xi_{i_k} \leq lQ < \sum\limits_{k=1}^{j}\xi_{i_k}\right)c_{0i_j}
\end{align*}
\endgroup
Since the expected recourse cost at customers $i_1$ to $i_{h-1}$ and at customers $i_{h+2}$ to $i_{|S|}$ is identical in both paths, the inequality $\mathcal{Q}(p) \geq \mathcal{Q}(p')$ simplifies to:
\begingroup
\allowdisplaybreaks
\begin{align*}
    & \mathcal{Q}(p) \geq \mathcal{Q}(p')\\
    \iff & 2\sum\limits_{l=1}^{\infty}\mathds{P}\left(\tilde{\xi} \leq lQ < \tilde{\xi} + \xi_{i_h}\right)c_{0i_h} + 2\sum\limits_{l=1}^{\infty}\mathds{P}\left(\tilde{\xi} + \xi_{i_h} \leq lQ < \tilde{\xi} + \xi_{i_h} + \xi_{i_{h+1}} \right)c_{0i_{h+1}} \geq\\
    & 2\sum\limits_{l=1}^{\infty}\mathds{P}\left(\tilde{\xi} \leq lQ < \tilde{\xi} + \xi_{i_{h+1}}\right)c_{0i_{h+1}} + 2\sum\limits_{l=1}^{\infty}\mathds{P}\left(\tilde{\xi} + \xi_{i_{h+1}} \leq lQ < \tilde{\xi} + \xi_{i_{h+1}} + \xi_{i_h} \right)c_{0i_h}\\
    \iff & \sum\limits_{l=1}^{\infty}\left( \mathds{P}\left(\tilde{\xi} \leq lQ < \tilde{\xi} + \xi_{i_h}\right) - \mathds{P}\left(\tilde{\xi} + \xi_{i_{h+1}} \leq lQ < \tilde{\xi} + \xi_{i_{h+1}} + \xi_{i_h} \right) \right)c_{0i_h} \geq\\
    & \sum\limits_{l=1}^{\infty} \left(\mathds{P}\left(\tilde{\xi} \leq lQ < \tilde{\xi} + \xi_{i_{h+1}}\right) - \mathds{P}\left(\tilde{\xi} + \xi_{i_h} \leq lQ < \tilde{\xi} + \xi_{i_h} + \xi_{i_{h+1}} \right)\right)c_{0i_{h+1}}\\
    \iff & \sum\limits_{l=1}^{\infty}\left( \mathds{P}\left(\tilde{\xi} + \xi_{i_h} > lQ\right) - \mathds{P}\left(\tilde{\xi} > lQ\right) - \mathds{P}\left(\tilde{\xi} + \xi_{i_{h+1}} + \xi_{i_h} > lQ\right) + \mathds{P}\left(\tilde{\xi} + \xi_{i_{h+1}} > lQ \right) \right)c_{0i_h} \geq\\
    & \sum\limits_{l=1}^{\infty} \left(\mathds{P}\left( \tilde{\xi} + \xi_{i_{h+1}} > lQ\right) - \mathds{P}\left(\tilde{\xi} > lQ \right) - \mathds{P}\left(\tilde{\xi} + \xi_{i_h} + \xi_{i_{h+1}} > lQ\right) + \mathds{P}\left(\tilde{\xi} + \xi_{i_h} > lQ \right)\right)c_{0i_{h+1}}\\
    \iff & \sum\limits_{l=1}^{\infty}\left( - \mathds{P}\left(\tilde{\xi} + \xi_{i_h} \leq lQ\right) + \mathds{P}\left(\tilde{\xi} \leq lQ\right) + \mathds{P}\left(\tilde{\xi} + \xi_{i_{h+1}} + \xi_{i_h} \leq lQ\right) - \mathds{P}\left(\tilde{\xi} + \xi_{i_{h+1}} \leq lQ \right) \right)c_{0i_h} \geq\\
    & \sum\limits_{l=1}^{\infty} \left(- \mathds{P}\left( \tilde{\xi} + \xi_{i_{h+1}} \leq lQ\right) + \mathds{P}\left(\tilde{\xi} \leq lQ \right) + \mathds{P}\left(\tilde{\xi} + \xi_{i_h} + \xi_{i_{h+1}} \leq lQ\right) - \mathds{P}\left(\tilde{\xi} + \xi_{i_h} \leq lQ \right)\right)c_{0i_{h+1}}\\
    \iff & (c_{0i_{h+1}} - c_{0i_{h}})\sum\limits_{l=1}^{\infty}\left( \mathds{P}\left(\tilde{\xi} + \xi_{i_h} \leq lQ\right) - \mathds{P}\left(\tilde{\xi} \leq lQ\right) - \mathds{P}\left(\tilde{\xi} + \xi_{i_{h+1}} + \xi_{i_h} \leq lQ\right) + \mathds{P}\left(\tilde{\xi} + \xi_{i_{h+1}} \leq lQ \right) \right) \geq 0\\
    \impliedby & \mathds{P}\left(\tilde{\xi} + \xi_{i_h} \leq lQ\right)  - \mathds{P}\left(\tilde{\xi} + \xi_{i_{h+1}} + \xi_{i_h} \leq lQ\right) \geq \mathds{P}\left(\tilde{\xi} \leq lQ\right) - \mathds{P}\left(\tilde{\xi} + \xi_{i_{h+1}} \leq lQ \right)  \hspace{2cm} , \forall l\geq 1
\end{align*}
\endgroup

Since $S$ respects the monotonicity condition, the last inequality, which corresponds to inequality (\ref{ineq:monotonicity_condition_1}) with $\tilde{S} = \{i_1,\dots,i_{h-1}\}$ and $(a,b)=(i_{h}, i_{h+1})$, is verified for any $l\in \mathds{N}$.

This implies that, starting from any path $p \in \mathcal{P}^1(S)$, one can iteratively perform inversions of consecutive customers $(a,b) \in S \times S$ such that $c_a < c_b$ until all the customers are sorted in non-increasing order of their distance to the depot without increasing the expected recourse cost. Therefore, the set $\argmin_{p \in \mathcal{P}^1(S)}\mathcal{Q}(p)$ contains at least one path $\tilde{p} = (\tilde{i}_1,\dots,\tilde{i}_{|S|})$ such that $c_{0\tilde{i}_1} \geq c_{0\tilde{i}_2} \geq \dots \geq c_{0\tilde{i}_{|S|}}$.

To conclude that $p^* \in \argmin_{p \in \mathcal{P}^1(S)}\mathcal{Q}(p)$, and thus that $\mathcal{Q}(p^*)$ is a valid lower bound on the expected recourse of any path $p \in \mathcal{P}^1(S)$, it remains to show that $\mathcal{Q}(\tilde{p}) = \mathcal{Q}(p^*)$. To do so, we sort the distances separating the customers of $S$ from the depot and remove the duplicates to obtain an ordered list $c_{0(1)} > \dots > c_{0(t)}$ of length $t \leq |S|$. For each $s \in \{1,\dots,t\}$, let us denote by $C_s = \{i \in S : c_{0i}=c_{0(s)}\}$ the set of customers that share the same distance $c_{0(s)}$ from the depot and by $\xi^s = \sum_{i \in C_s} \xi_i$ the sum of their demands. By construction, in both paths $\tilde{p}$ and $p^*$, a permutation of $C_1$ will first be visited, followed by a permutation of $C_2$ and so on until $C_t$. Denoting the number of customers in group $C_s$ as $v_s = |C_s|$ and the total number of customers visited before moving to group $C_{s+1}$ as $V_s = \sum_{q=1}^s v_r$, we can demonstrate the expected equality.
\begingroup
\allowdisplaybreaks
\begin{align}
    \label{lower_bound_sorted_e1} \mathcal{Q}(\tilde{p}) &= 2\sum_{j=1}^{|S|}\sum_{l=1}^{\infty}\mathds{P}\left(\sum_{k=1}^{j-1}\xi_{\tilde{i}_k} \leq lQ < \sum_{k=1}^{j}\xi_{\tilde{i}_k}\right)c_{0\tilde{i}_j} \\
    \label{lower_bound_sorted_e2} &= 2\sum_{s=1}^{t}\sum_{j = V_{s-1}+1}^{V_s}\sum_{l=1}^{\infty}\mathds{P}\left(\sum_{k=1}^{j-1}\xi_{\tilde{i}_k} \leq lQ < \sum_{k=1}^{j}\xi_{\tilde{i}_k}\right)c_{0(s)} \\
    \label{lower_bound_sorted_e3} &= 2\sum_{s=1}^{t}c_{0(s)}\sum_{l=1}^{\infty}\left(\sum_{j = V_{s-1}+1}^{V_s}\mathds{P}\left(\sum_{k=1}^{j-1}\xi_{\tilde{i}_k} \leq lQ < \sum_{k=1}^{j}\xi_{\tilde{i}_k}\right)\right) \\
    \label{lower_bound_sorted_e4} &= 2\sum_{s=1}^{t}c_{0(s)}\sum_{l=1}^{\infty}\mathds{P}\left(\sum_{q=1}^{s-1}\xi^{q} \leq lQ < \sum_{q=1}^{s}\xi^{q}\right) \\
    \label{lower_bound_sorted_e5} &= 2\sum_{s=1}^{t}c_{0(s)}\sum_{l=1}^{\infty}\left(\sum_{j = V_{s-1}+1}^{V_s}\mathds{P}\left(\sum_{k=1}^{j-1}\xi_{i^*_k} \leq lQ < \sum_{k=1}^{j}\xi_{i^*_k}\right)\right) \\
    \label{lower_bound_sorted_e6} &= 2\sum_{s=1}^{t}\sum_{j = V_{s-1}+1}^{V_s}\sum_{l=1}^{\infty}\mathds{P}\left(\sum_{k=1}^{j-1}\xi_{i^*_k} \leq lQ < \sum_{k=1}^{j}\xi_{i^*_k}\right)c_{0(s)} \\
    \label{lower_bound_sorted_e7} &= 2\sum_{j=1}^{|S|}\sum_{l=1}^{\infty}\mathds{P}\left(\sum_{k=1}^{j-1}\xi_{i^*_k} \leq lQ < \sum_{k=1}^{j}\xi_{i^*_k}\right)c_{0i^*_j} \\
    \label{lower_bound_sorted_e8} &= \mathcal{Q}(p^*)
\end{align}
\endgroup

The equivalence of equalities (\ref{lower_bound_sorted_e3}) and (\ref{lower_bound_sorted_e4}) and, analogously, of equalities (\ref{lower_bound_sorted_e4}) and (\ref{lower_bound_sorted_e5}), comes from the fact that the $l$-th restocking trip will occur at a customer of set $C_s$ if and only if the total demand of the customers of sets $C_1$ to $C_{s-1}$ is less than or equal to $lQ$ and the total demand of the customers of sets $C_1$ to $C_{s}$ exceeds $lQ$. The sum of the probabilities of the mutually exclusive events that are considered for each $s \in \{1,\dots,t\}$ and for each $l \in \mathds{N}$ in both equations (\ref{lower_bound_sorted_e3}) and (\ref{lower_bound_sorted_e5}) is thus equal to the probability of their disjunction, which is given in (\ref{lower_bound_sorted_e4}).
\endproof

The definition of the single-route lower bound follows.
\begin{equation}
    L_1(S)=\mathcal{Q}(p^*)\label{eq:Single-route lb}
\end{equation}

\subsubsection{Dynamic programming algorithm.} \label{sec:DP for Recourse}

The lower bound proposed in this section applies to sets of customers $S$ that can be assigned to $m\geq 2$ vehicles. It is inspired by the lower bound of \cite{laporte2002integer}, where the demands are seen as continuous quantities that must be distributed among the vehicles. This new bound, which requires the monotonicity of the recourse function, improves over \cite{laporte2002integer} in two ways. First, the value of the lower bound can be greatly increased by considering that multiple failures can happen for a vehicle. This is done by computing for each vehicle $k$ and for each quantity of expected demand $q \in [0, Q]$ the least-cost sequence of customers by taking advantage of Proposition \ref{prop:lower_bound}. The second improvement is that our search for an assignment of the demands can be restricted to integer solutions through the use of a new dynamic programming algorithm. The new bound is computed in two steps.

In the first step, the least-cost sequence for assigning $q_h \in [0, Q]$ units of expected demand to a vehicle using the $|S|-h+1$ farthest customers of $S$ is computed. To do so, the customers of $S$ are sorted by non-decreasing distance to the depot. The least-cost sequence is computed by dynamic programming over a set of $|S|$ stages indexed by $h= 1$ to $|S|$. The state $q_h$ represents the number of units of expected demand assigned in stage $h$. In state $q_h$, we denote the total demand of the assigned customers by the RV $\xi_{q_h}$. Its distribution is given by $\xi_{q_h} \sim \text{Poisson}(q_h)$ if the demands are modeled as Poisson variables and by $\xi_{q_h} \sim \mathcal{N}(q_h, var(q_h))$ if the demands are normally distributed. In the case of normal demands, we set $var(q_h) = \min\{\sum_{i=h}^{|S|}\sigma^2_i z_i|\sum_{i=h}^{|S|}\mu_iz_i=q_h, z_i \in \{0,1\}\}$, which is the smallest possible variance for the total demand of the assigned customers. Function $G_h(q_h)$ finds the value of the least-cost sequence.
\begin{equation*}
    G_h(q_h) = \begin{cases}
        \min\left \{G_{h+1}(q_h), G_{h+1}(q_h-\mu_h) + 2 c_{0h} \sum\limits_{l=1}^\infty \mathds{P}(\xi_{q_h-\mu_h} \leq lQ < \xi_{q_h} ) \right\}, \label{dp_l2} & \mbox{for } h < |S| \mbox{ and } q_h \geq \mu_h, \\
        G_{h+1}(q_h), & \mbox{for } h < |S| \mbox{ and } q_h < \mu_h, \\
        0, & \mbox{for } h = |S| \mbox{ and } q_h = 0, \\
        2 c_{0h} \mathds{P}(\xi_{q_h} > Q), & \mbox{for } h = |S| \mbox{ and } q_h = \mu_h, \\
        \infty, & \mbox{for } h = |S| \mbox{ and } q_h \notin \{0, \mu_h\}.
        \end{cases}
\end{equation*}

The first case selects the lowest-cost assignment between the previous stage in state $q_h$ and the previous stage in state $q_h - \mu_h$ plus the expected recourse cost paid at customer $h$. The second case indicates that customer $h$ cannot be used to fulfill $q_h < \mu_h$ units of expected demand, only those after $h$ can thus be used. The last three cases cover the initial stage. The farthest customer $h=|S|$ can be either visited or not, which corresponds to assigning either 0 or $\mu_{h}$ units of expected demand to the vehicle.

In the second step, the bound is obtained by computing the least-cost assignment of the $\sum_{i \in S} \mu_i$ units of expected demand to the vehicles. This is done by defining another dynamic program that uses the previously computed functions $G_h(q)$ and only allows the $h$-th customer to be assigned to vehicle $k$ if $k \leq h$. This dynamic program is defined over $m$ stages indexed by $k=1$ to $m$. The state $q_k$ represents the number of units of expected demand that are assigned to the first $k$ vehicles. At stage $k$, variable $x_k$ decides which quantity is assigned to vehicle $k$. Then, function $H_k(q_k)$ finds the optimal assignment of $q_k$ units to the first $k$ vehicles. It is defined as follows.
\begin{equation*}
    H_k(q_k)= \begin{cases}
      \infty, & \mbox{if } q_k > kQ, \\
      G_k(q_k), & \mbox{for } k = 1 \mbox{ and } q_k \leq Q, \\
  		\min\limits_{0\leq x_k \leq \min(q_k,Q)} \{ G_k(x_k)+H_{k-1}(q_k-x_k)\}, & \mbox{if } 2 \leq k \leq m.
  	\end{cases}
\end{equation*}

Finally, the expected recourse cost of the optimal assignment of the demands is given by:
\begin{equation}
  L_2(S) = H_{m} \left(D\right), \label{eq:L2}
\end{equation}
where $D = \sum_{i \in S}\mu_i$. This value can be computed in $O(m D Q)$.

\subsubsection{A set-covering formulation.}\label{sec:spf}

The last lower bound is based on a set-covering formulation, where the linear relaxation is solved by column generation. Each column represents a path that visits some customers of $S$ and respects the expected capacity constraint. The problem is to find a set of feasible paths of cardinality $m$ that visits all customers of $S$ and minimizes the total expected recourse cost. Let $\mathcal{P}(S)$ be the set of all feasible paths visiting only customers of $S$. We define $z_p$ as a binary variable indicating whether path $p \in \mathcal{P}(S)$ is taken or not. Also, the binary coefficient $b_{jp}$ indicates whether customer $j$ is in path $p$. The formulation is as follows:
\begingroup
\allowdisplaybreaks
\begin{align}
   L_3(S) = \min\ & \sum_{p \in \mathcal{P}(S)} \mathcal{Q}(p) z_p  \label{SetCoveringCG1}\\
    \text{s.t.}\ & \sum_{p \in \mathcal{P}(S)} b_{ip}z_p \geq 1, & i \in S,\label{SetCoveringCG2}\\
    & \sum_{p \in \mathcal{P}(S)}z_p = m,\label{SetCoveringCG3}\\
    & z_p \in \{0,1\},&  p \in \mathcal{P}(S). \label{SetCoveringCG4}
\end{align}
\endgroup

The objective function \eqref{SetCoveringCG1} minimizes the sum of the expected recourse costs of the selected paths. Constraints \eqref{SetCoveringCG2} state that each customer of $S$ must be visited. Constraint \eqref{SetCoveringCG3} ensures that exactly $m$ paths are selected, and constraints \eqref{SetCoveringCG4} define the domain of the variables.

As the set of paths $\mathcal{P}(S)$ is too large to be generated up front, we rely on column generation to solve the problem. An initial set of columns is first generated to obtain a feasible solution for the linear relaxation. Then, a pricing problem is solved to find a column having a negative reduced cost. If such a column is found, it is added to the model, and the linear relaxation is solved again. This process continues until no more negative reduced cost columns can be found.

The pricing problem of finding the least-cost path can be formulated using the dynamic program of the previous section. The first equation, when adding customer $h$, is modified by subtracting the value of the dual variable associated with the constraint \eqref{SetCoveringCG2} of customer $h$. A negative cost column is found if the value of the least-cost sequence is less than the value of the dual variable associated with constraint \eqref{SetCoveringCG3}.

\section{Implementation of the DL-shaped method}

This section details the steps performed when executing the DL-shaped method.

\begin{enumerate}
    \item Compute $L_m = L_3(N)$ for each number of vehicles $m\in M$.
    \item Build model (\ref{model:depot_degree})-(\ref{model:int}) with the objective function (\ref{model1:vrp7}) and add constraint \eqref{eq:L}.
    \item Solve the linear relaxation of the model.
    \item Verify if there are any violated subtour or capacity constraints \eqref{model:cap} using the CVRPSEP package of \cite{lysgaard2003cvrpsep}. Add the violated inequalities \eqref{model:cap} to the model. For each identified set $S$, check if its associated S-Cut is violated, and if so, add the inequality to the model.
    \item If the solution is fractional, enumerate its connected components. For each component, verify if an S-Cut is violated, and if one is found, add the S-Cut to the model. If the number of nodes in the component equals the number of non-zero valued edges plus one, then the component is a path. Then, check if the P-Cut associated with the path is violated, and if so, add it to the model.
    \item Go back to Step 3 if violated inequalities were found in Steps 4 or 5. Otherwise, branch on a fractional edge variable if the solution is fractional and return to Step 3.
    \item For each route in the solution, add the violated P\&S-Cuts to the model. Return to Step 3 if violated inequalities were found.
\end{enumerate}

The procedure stops when all the nodes of the brand-and-bound tree have been explored. The incumbent solution is then optimal. In our implementation, the S-Cut associated with a set of customers $S$ is generated using the lower bound $L_1(S)$ if the sum of their expected demands does not exceed the vehicle capacity. Otherwise, the lower bound $L_2(S)$ is used, with $m=\lceil (\sum_{i \in S}\mu_i)/Q \rceil$. To compute the P\&S-Cuts of Step 7, we iterate on the routes of the current solution and generate all the subpaths in which at most five customers are removed from the original route. The recourse cost and the single-route lower bound $L_1$ are then computed for each subpath. The five most violated P-Cuts and the five most violated S-Cuts are added to the model.

\section{Computational experiments}\label{sec:Numerical results}

This section analyzes the performance of the DL-shaped method through computational experiments. All the results were produced with a C\texttt{++} implementation with Cplex 12.10 on a machine with an AMD Rome 7532 @ 2.40GHz CPU. For all the computational experiments, we defined a maximum time of one hour per instance. A quick upper bound is also provided to the solver by applying the adaptive large neighborhood search of \cite{Pisinger2006}. The detailed results can be found at \href{https://sites.google.com/view/jfcote/}.

Section \ref{sec:Description of the instances} describes the existing instances that are used in our experiments and introduces a new set of instances. In Section \ref{sec:results_lbs}, the new lower bounds on the recourse are compared with classical bounds from the literature. Sections \ref{sec:results_literature} compare our algorithm to previous methods on two different instance sets. Finally, Section \ref{sec:Results of New Instances} presents the performance of our method on the new set of instances and summarizes the key takeaways of our experiments.

\subsection{Characterization of the instances}\label{sec:Description of the instances}

Computational experiments are based on three sets of instances: two from the literature and a new one. The first set is taken from \cite{jabali2014partial} and contains 270 instances. In each instance, $n$ ranges from 40 to 80, and the demands are normally distributed, with parameters $(\mu_i,\sigma_i)$ such that $\mu_i=3\sigma_i$ and $\mu_i \in \{1,\dots,10\}$. The number of vehicles is fixed to $\bar{m} \in \{2,3,4\}$ and the filling coefficient $\bar f = \frac{\sum_{i \in N} \mu_i}{\bar{m}Q}$ ranges from 0.80 to 0.95. The DL-shaped method can be used on this set since they respect the monotonicity property. We verified that inequality (\ref{ineq:monotonicity_condition_2}) holds for each capacity $Q$ in this set of instances, for each number of restocking trips $l \in \{1,2,3\}$, and for each possible distribution of $\xi_a \sim \mathcal{N}(\mu_a, \sigma^2_a)$, $\xi_b \sim \mathcal{N}(\mu_b, \sigma^2_b)$, and $\tilde{\xi} \sim \mathcal{N}(\tilde{\mu},\tilde{\sigma}^2)$. We ignored the case $l \geq 4$ since the probability of more than three restocking trips to occur is negligible.

The second set is taken from classical CVRP instances and contains the 95 instances of the sets A, B, E, F, M, and P of the \texttt{CVRPLIB} repository \citep{uchoa2017new}. In these instances, the fleets are unlimited, and the customer demands are modeled as Poisson RVs. The expectation $\mu_i$ of each customer $i \in N$ is set to the deterministic demand of the CVRP instance. Also, as done by \cite{christiansen2007branch}, \cite{gauvin2014branch} and \cite{florio2020new}, demands and the capacity of each instance are divided by their greatest common divisor.

For the new set, we generated 1980 instances using the method of \cite{jabali2014partial}. To obtain a challenging set of instances, we produced ten instances for each combination of the following problem parameters: $n\in\{20,30,40,50,60,70,80,90,100,110,120\}$, $\bar{m}\in\{2,3,4,5,6,7\}$, and $\bar{f}\in\{0.85,0.90,0.95\}$. The coefficient of dispersion of the demand of each customer $i \in N$ is given by $D=\sigma^2_i/\mu_i=1$. By Proposition \ref{prop:normal}, these instances therefore respect the monotonicity property.

\subsection{Lower bounds on the recourse}\label{sec:results_lbs}

This section compares our new lower bounds on the recourse to those previously used in the literature based on the relative gap they provide when applied to an optimal first-stage solution. This measure is defined as follows.
\begin{equation*}
    \text{Gap\%}=\frac{\text{Recourse in the best or optimal solution}-\text{lower bound}}{\text{Recourse in the best or optimal Solution}}\times 100
\end{equation*}

For the first set of instances of \cite{jabali2014partial}, the results are grouped by number of customers, whereas the results are presented separately for each set A, B, E, M, and P for the CVRP instances. In Tables \ref{tab:bounds_jabali} and \ref{tab:bounds_cvrp}, the columns LLVH02 and LSG18 present the gaps obtained from the implementation of the lower bound from \cite{laporte2002integer} and \cite{louveaux2018exact}, respectively. The performance of our three lower bounds, $L_1$, $L_2$, and $L_3$ is also reported. Column $L_1$ presents the gaps obtained from computing $\sum_{r \in \mathcal{R}^*}L_1(r)$, where $\mathcal{R}^*$ is the set of routes obtained in either the optimal or best solution found at the end of the DL-shaped method. Columns $L_2$ and $L_3$ are obtained by computing $L_2(N)$ and $L_3(N)$. The average computing time of $L_3(N)$ is also reported. The computing time of all the other bounds is negligible.

\begin{table}[ht]
\centering
\begin{tabular}{|c c|c c c c c c|}
\hline
$n$ & \# & LLVH02 & LSG18 & $L_1$ & $L_2$ & $L_3$ & $L_3$ (s) \\ \hline
40 & 30 & 99.9 & 99.9 & 3.5 & 89.7 & 60.0 & 0.9 \\
50 & 60 & 96.0 & 99.6 & 6.9 & 85.1 & 55.1 & 2.9 \\
60 & 90 & 88.6 & 97.2 & 7.4 & 78.5 & 50.8 & 6.4 \\
70 & 60 & 86.3 & 96.7 & 9.0 & 73.8 & 44.2 & 13.7 \\
80 & 30 & 76.6 & 94.5 & 10.2 & 61.4 & 30.4 & 21.7 \\ \hline
\multicolumn{2}{|c|}{Average} & 89.5 & 97.6 & 7.4 & 77.7 & 48.1 & 8.8 \\ \hline
\end{tabular}
\caption{Relative gap of lower bounds on the instances of \cite{jabali2014partial}}
\label{tab:bounds_jabali}
\end{table}

\begin{table}[ht]
\centering
\begin{tabular}{|c c|c c c c c c|} \hline
Set & \# & LLVH02 & LSG18 & $L_1$ & $L_2$ & $L_3$ & $L_3$ (s) \\ \hline
A & 27 & 67.6 & 100.0 & 1.9 & 67.0 & 64.9 & 0.5 \\
B & 23 & 58.4 & 99.9 & 0.5 & 59.6 & 55.9 & 0.5 \\
E & 13 & 60.4 & 99.5 & 4.0 & 59.4 & 56.9 & 0.7 \\
M & 3 & 97.6 & 100.0 & 3.2 & 100.0 & 93.8 & 6.3 \\
P & 5 & 50.8 & 99.8 & 2.2 & 49.3 & 43.3 & 7.0 \\
F & 24 & 38.5 & 95.5 & 2.0 & 37.2 & 35.2 & 0.4 \\ \hline
\multicolumn{2}{|c|}{Average} & 57.1 & 98.8 & 1.9 & 56.8 & 53.9 & 1.0 \\ \hline
\end{tabular}
\caption{Relative gap of lower bounds on the CVRP instances}
\label{tab:bounds_cvrp}
\end{table}

The lower bound of \cite{louveaux2018exact} is systematically the weakest, followed by that of \cite{laporte2002integer}. Their average gap increase with both the fleet size and the number of customers. The bound of \cite{louveaux2018exact} assumes a single large vehicle and therefore becomes a poor approximation of the recourse as the fleet size increases. Regarding the bound of \cite{laporte2002integer}, it suffers from an underestimation of the variance of the total demand on a route, which worsens as the number of customers grows.

These results demonstrate that our new bounds significantly improve over the previous bounds from the literature. The set covering bound $L_3$ achieves an average gap of approximately $50\%$ on average on both sets of instances. The DP bound $L_2$ is generally slightly less tight than $L_3$ but is significantly cheaper to compute. The good computational performance of $L_2$ can be attributed to the fact that the least-cost sequence can be computed by sorting the customers according to their distance to the depot. The single-route lower bound $L_1$ is by far the tightest, with an average gap of $7.4\%$ for the first set and $1.4\%$ for the second. For a given set of customers, the quality of bounds $L_1$ and $L_2$ directly depends on the similarity of the least-cost sequence for the recourse cost provided by Proposition \ref{prop:lower_bound} and the total least-cost sequence that takes the first-stage routing cost into account. In particular, the performance of bound $L_1$ reported in Table \ref{tab:bounds_cvrp} indicates the substantial similarity of these sequences for the CVRP instances.

\subsection{Results for the instances of the literature }\label{sec:results_literature}

This section reports results on the two sets of instances of the literature. It compares the efficiency of the DL-shaped method against other B\&C methods and against B\&P methods.

Table \ref{tab:Comparison jabali} presents the results against the integer L-shaped methods of \cite{jabali2014partial} and \cite{Ymro2023}. The performance of the three methods is presented under the columns JRGL14, HS22, and DL-Shaped of Table \ref{tab:Comparison jabali}. Each row of Table \ref{tab:Comparison jabali} contains 30 instances with the same number of customers $n$ and the same fleet size $\bar{m}$. The number of optimal solutions found by each method and their average computing time are reported, as well as the average optimality gap of the unsolved instances.

\begin{table}[H]
    \centering
    \begin{tabular}{|c c|c c c|c c c|c c|}
    \hline
 &  & \multicolumn{3}{c|}{JRGL14} & \multicolumn{3}{c|}{HS22} & \multicolumn{2}{c|}{DL-Shaped} \\\hline
n & $\bar{m}$ & Opt & Sec & Gap\% & Opt & Sec & Gap\% & Opt & Sec \\\hline
40 & 4 & 9 & 1240.7 & 1.5 & 28 & 91.9 & 1.6 & 30 & 2.0 \\
50 & 3 & 16 & 6918.0 & 0.7 & 29 & 101.5 & 2.9 & 30 & 9.8 \\
50 & 4 & 5 & 1360.8 & 1.9 & 25 & 224.8 & 1.7 & 30 & 30.9 \\
60 & 2 & 24 & 1393.0 & 0.4 & 30 & 72.2 & - & 30 & 1.2 \\
60 & 3 & 6 & 2766.0 & 0.7 & 27 & 191.4 & 1.2 & 30 & 16.7 \\
60 & 4 & 3 & 4922.0 & 2.0 & 25 & 262.1 & 1.9 & 30 & 25.3 \\
70 & 2 & 17 & 2577.5 & 0.5 & 30 & 99.2 & - & 30 & 2.0 \\
70 & 3 & 9 & 1753.3 & 1.5 & 24 & 563.7 & 1.5 & 30 & 46.2 \\
80 & 2 & 13 & 1809.2 & 0.5 & 28 & 193.7 & 1.0 & 30 & 19.8 \\ \hline
\multicolumn{2}{|c|}{Total or avg.} & 102 & 2711.4 & 1.2 & 246 & 185.7 & 1.6 & 270 & 17.1 \\ \hline
    \end{tabular}
    \caption{Performance comparison with \cite{jabali2014partial, Ymro2023} on the instances of \cite{jabali2014partial}}
     \label{tab:Comparison jabali}
\end{table}

Table \ref{tab:Comparison jabali} indicates that the DL-shaped method achieves state-of-the-art results by solving to optimality all the 270 instances of the first set, while  \cite{jabali2014partial} and \cite{Ymro2023} respectively solve 102 and 246 instances. Furthermore, our average resolution time is 17.1 seconds per instance, compared to 2711.4 and 185.72 seconds for the instances that were solved to optimality by \cite{jabali2014partial} and \cite{Ymro2023}, respectively. These results highlight the advantage of disaggregating the recourse in the context of the integer L-shaped method, as well as the capacity of our LBFs to approximate the recourse function efficiently.

In Table \ref{tab:Comparison gauvin}, the DL-shaped method is compared to the algorithms of \cite{christiansen2007branch} (CL07), \cite{gauvin2014branch} (GDG14), \cite{florio2020new} (FHM20), and \cite{Ymro2023} (HS23) on the second set of instances.
\begin{table}[H]
\centering
\begin{tabular}{|c|c|c c c c |c c c|}
\hline
Set & \# & CL07 & GDG14 & FHM20* & HS23* & DL-Shaped & Gap\% & Time (s) \\ \hline
A & 27 & 6/19 & 22/26 & 15/27 & 1 & 9 & 6.8 & 202.1 \\
B & 23 & - & 7/23 & - & - & 6 & 4.3 & 1067.1 \\
E & 13 & 2/3 & 5/11 & 3/8 & 6 & 7 & 7.3 & 127.2 \\
F & 3 & - & - & - & - & 2 & 1.8 & 18.7 \\
M & 5 & - & 1/1 & - & - & 0 & 9.9 &- \\
P & 24 & 11/18 & 20/22 & 17/23 & 7 & 11 & 6.0 & 143.9 \\\hline
Total or avg. & 95 & 19/40 & 55/83 & 35/58 & 14 & 35 & 6.1 & 306.7 \\\hline
Max. & $n/m$ &  6.3 & 18.5 & 10.5 & 10.5 & 25 & & \\\hline
\end{tabular}
\caption{Performance comparison with \cite{christiansen2007branch, gauvin2014branch, florio2020new} on the CVRP instances. * indicate that the authors use the OR policy.}
\label{tab:Comparison gauvin}
\end{table}

A rigorous comparison of these algorithms is difficult to achieve, as the previous works do not all report results for the complete set of instances, and some use the OR policy instead of the DTD policy. Nevertheless, the fact that the DL-shaped method can solve more than a third of the 95 instances of this set is a new milestone for a B\&C approach. Indeed, the previous integer L-shaped algorithms from the literature have not been successful in solving these instances. \cite{jabali2014partial, louveaux2018exact, salavati2019exact} did not report any results, and HS23 can only solve 14 of them under the OR policy. Our method achieves comparable performance as GDG14 for the B and E sets and generally produces solutions of good quality, with the average gaps being at most 10\% for each set. Also, the DL-shaped method can solve two of the three instances of set F. In one of those, the vehicle capacity is $30,000$. For comparison, the highest vehicle capacity is $8,000$ among the instances solved to optimality by HS23, and $3,000$ for FHM20. These results illustrate the advantage of the DL-shaped method over previous algorithms from the literature for instances with very high vehicle capacities and, more generally, that B\&C methods are less sensitive than B\&P methods to large vehicle capacities and expected customer demands. The last aspect in which our method stands out is the length of the routes that can be handled. The last row of Table \ref{tab:Comparison gauvin} reports the maximal ratio of customers to vehicles among the instances that can be solved to optimality by each method. This value is 25 for the DL-shaped method, more than twice that of FHM20 and HS23. It was achieved on an instance of set P that contains 100 customers and whose optimal solution comprises four routes, the longest of which visits 29 customers.

\subsection{Results for the new instances }\label{sec:Results of New Instances}

This section presents the results for our new set of instances. For each number of customers and vehicles, Table \ref{tab:New instances} reports the number of instances, out of 30, that could be solved within the time limit and their average resolution time. In total, 1175 instances were solved, and the average optimality gap of the 805 unsolved instances was less than $5\%$ in each group. Even for as many as 120 customers, instances with two vehicles were generally solved easily, with an average time of one minute per solved instance. Although the difficulty of the new instances appears to increase with the number of vehicles, some of the instances that were solved by the DL-shaped method are significantly larger than those previously found in the literature. These include two instances with 50 customers and seven vehicles and one instance with 120 customers and five vehicles.

%


\begin{table}[H]
    \begin{center}
    \makebox[0cm]{
    \begin{tabular}{|c|c c c c c c|c|c c c c c c|c|}
    \hline
&\multicolumn{7}{c|}{Number of optimal solutions}&\multicolumn{7}{c|}{Average computing times (s)}\\ \hline
$n \backslash \bar{m}$ & 2 & 3 & 4 & 5 & 6 & 7 & All & 2 & 3 & 4 & 5 & 6 & 7 & All \\ \hline
20 & 30 & 30 & 30 & 30 & 30 & 30 & 180 & 0.0 & 0.4 & 1.3 & 7.9 & 32.4 & 10.7 & 8.8 \\
30 & 30 & 30 & 30 & 30 & 24 & 20 & 164 & 0.1 & 0.6 & 5.6 & 276.1 & 807.3 & 863.9 & 275.2\\
40 & 30 & 30 & 30 & 28 & 13 & 3 & 134 & 0.3 & 1.4 & 42.2 & 498.0 & 867.0 & 1718.3 & 236.5\\
50 & 30 & 30 & 28 & 23 & 4 & 2 & 117 & 0.6 & 45.1 & 232.2 & 516.1 & 913.6 & 1297.1 & 222.1\\
60 & 30 & 30 & 24 & 17 & 3 & 0 & 104 & 3.9 & 12.3 & 275.3 & 1114.8 & 1926.9 & - & 306.0 \\
70 & 30 & 30 & 22 & 7 & 1 & 0 & 90 & 1.6 & 286.1 & 272.3 & 1653.7 & 400.0 & - & 295.5 \\
80 & 30 & 30 & 19 & 10 & 2 & 0 & 91 & 17.7 & 198.1 & 1024.2 & 1455.7 & 1667.5 & - & 481.6 \\
90 & 30 & 28 & 19 & 6 & 1 & 0 & 84 & 15.8 & 233.3 & 710.4 & 611.6 & 2184.9 & - &  313.8\\
100 & 30 & 25 & 18 & 6 & 0 & 0 & 79 & 16.3 & 396.9 & 956.0 & 1382.5 & -  & - & 454.6 \\
110 & 30 & 26 & 9 & 5 & 0 & 0 & 70 & 101.9 & 520.4 & 234.7 & 1484.3 & -  & - & 373.2 \\
120 & 29 & 21 & 11 & 1 & 0 & 0 & 62 & 62.2 & 306.3 & 317.4 & 1468.2 & -  & - & 212.8 \\ \hline
All & 329 &310 & 240 & 163 & 78 & 55 & 1175 & 19.9 & 170.1 & 318.1 & 615.2 & 602.2 & 460.9 & 262.3 \\ \hline
    \end{tabular}
    }
    \end{center}
    \caption{Solved instances and computation times of the DL-shaped method for the new set}
    \label{tab:New instances}
\end{table}

\begin{figure}[H]
\centering
\begin{tikzpicture}[scale=0.88]
	\begin{axis}[
		title={Largest solved instances},
        name=ax1,
		ylabel={Number of vehicles},
        y label style={at={(0.07,0.5)}},
		xlabel={Number of customers},
		ymin=0, ymax=8,
		xmin=0, xmax=130,
            ytick={1,2,3,4,5,6,7},
            xtick={20,40,60,80,100,120},
        legend style={font=\footnotesize, at={(-0.144,-0.23)},anchor=north west,legend columns=3, column sep=0.2cm, legend cell align=left},
		ymajorgrids=true,
            xmajorgrids=true,
		grid style=dashed,
		]
\addplot[ 
    black,
    mark=*,
    mark options={fill=black},
    visualization depends on=\thisrow{alignment} \as \alignment,
    nodes near coords, 
    point meta=explicit symbolic, 
    every node near coord/.style={anchor=\alignment} 
    ] table [
     meta index=2 
     ] {
x   y   label alignment
16	2	1.0   90
70	2	0.6   30
};
\addplot[ 
    black,
    thick,
    mark=oplus,
    mark options={fill=blue},
    visualization depends on=\thisrow{alignment} \as \alignment,
    nodes near coords, 
    point meta=explicit symbolic, 
    every node near coord/.style={anchor=\alignment} 
    ] table [
     meta index=2 
     ] {
x   y   label alignment
90	1	1.0    90
};
\addplot[
    black,
    mark=triangle,
    mark options={fill=white},
    visualization depends on=\thisrow{alignment} \as \alignment,
    nodes near coords, 
    point meta=explicit symbolic, 
    every node near coord/.style={anchor=\alignment} 
    ] table [
     meta index=2 
     ] {
x   y   label alignment
25	4	0.9	-90
50	3	0.9	-90
100	2	0.9	-90
};


\addplot[ 
    black,
    mark=square,
    visualization depends on=\thisrow{alignment} \as \alignment,
    nodes near coords, 
    point meta=explicit symbolic, 
    every node near coord/.style={anchor=\alignment} 
    ] table [
     meta index=2 
     ] {
x   y   label alignment
101	3	0.9	-90
};

\addplot[ 
     black,
     mark=o,
     visualization depends on=\thisrow{alignment} \as \alignment,
     nodes near coords, 
     point meta=explicit symbolic, 
     every node near coord/.style={anchor=\alignment} 
     ] table [
      meta index=2 
      ] {
x   y   label alignment
 60	4	0.85	-90
 70	3	0.9	-130
 80	2	0.95   120
 };

\addplot[ 
    black,
    mark=pentagon,
    visualization depends on=\thisrow{alignment} \as \alignment,
    nodes near coords, 
    point meta=explicit symbolic, 
    every node near coord/.style={anchor=\alignment} 
    ] table [
     meta index=2 
     ] {
x   y   label alignment
20	7	0.95	-90
50	7	0.9	-90
60	6	0.9	-120
90	6	0.85	-90
100	5	0.9	-120
120	5	0.85	-90
120	2	0.95	90
};

\legend{\cite{gendreau1995exact} (D),
        \cite{hjorring1999new} (D),
        \cite{laporte2002integer} (P),
        \cite{louveaux2018exact} (D),
        \cite{Ymro2023} (N),
        DL-shaped (N)}
\end{axis}

\begin{axis}[
		title={Smallest unsolved instances},
        at={(ax1.south east)},
        xshift=2cm,
		ylabel={Number of vehicles},
        y label style={at={(0.07,0.5)}},
		xlabel={Number of customers},
		ymin=0, ymax=8,
		xmin=0, xmax=130,
            ytick={1,2,3,4,5,6,7},
            xtick={20,40,60,80,100,120},
		ymajorgrids=true,
            xmajorgrids=true,
		grid style=dashed,
		]
\addplot[ 
    black,
    mark=*,
    mark options={fill=black},
    visualization depends on=\thisrow{alignment} \as \alignment,
    nodes near coords, 
    point meta=explicit symbolic, 
    every node near coord/.style={anchor=\alignment} 
    ] table [
     meta index=2 
     ] {
x   y   label alignment
16	2	0.6	-90
};
\addplot[ 
    black,
    thick,
    mark=oplus,
    mark options={fill=blue},
    visualization depends on=\thisrow{alignment} \as \alignment,
    nodes near coords, 
    point meta=explicit symbolic, 
    every node near coord/.style={anchor=\alignment} 
    ] table [
     meta index=2 
     ] {
x   y   label alignment
90	1	1.0    90
};
\addplot[
    black,
    mark=triangle,
    mark options={fill=white},
    visualization depends on=\thisrow{alignment} \as \alignment,
    nodes near coords, 
    point meta=explicit symbolic, 
    every node near coord/.style={anchor=\alignment} 
    ] table [
     meta index=2 
     ] {
x   y   label alignment
100	2	0.9	90
};


\addplot[ 
    black,
    mark=square,
    visualization depends on=\thisrow{alignment} \as \alignment,
    nodes near coords, 
    point meta=explicit symbolic, 
    every node near coord/.style={anchor=\alignment} 
    ] table [
     meta index=2 
     ] {
x   y   label alignment
31  3 0.85 -90
34  2 0.95 90
};
\addplot[ 
     black,
     mark=o,
     visualization depends on=\thisrow{alignment} \as \alignment,
     nodes near coords, 
     point meta=explicit symbolic, 
     every node near coord/.style={anchor=\alignment} 
     ] table [
      meta index=2 
      ] {
 x   y   label alignment
 40	4	0.85	-50
 50	3	0.9	90
 80	2	0.95   -130
 };
\addplot[ 
    black,
    mark=pentagon,
    visualization depends on=\thisrow{alignment} \as \alignment,
    nodes near coords, 
    point meta=explicit symbolic, 
    every node near coord/.style={anchor=\alignment} 
    ] table [
     meta index=2 
     ] {
x   y   label alignment
30	6	0.9	-90
40	5	0.85	-130
50	4	0.95	-130
80	4	0.85	-90
90	3	0.9	-120
100	3	0.9	-150
120	2	0.95	-120

};

\end{axis}
\end{tikzpicture}
\caption{Largest solved and smallest unsolved instances in the literature with their filling coefficient}
\label{fig1:ResFromLit}
\end{figure}
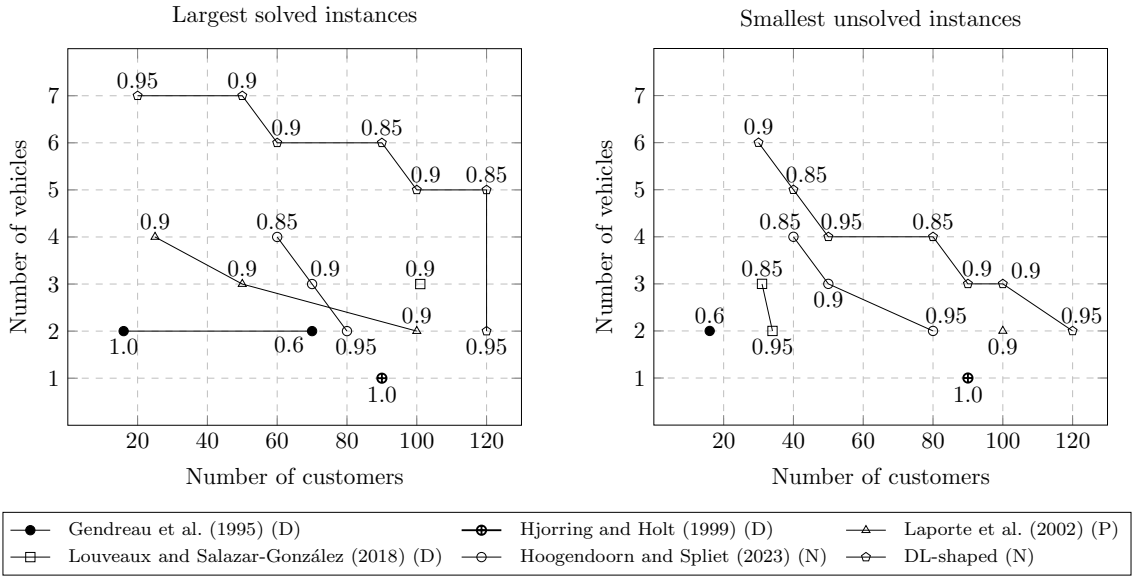

Figure \ref{fig1:ResFromLit} summarizes the parameters of the largest instances that could be solved and the smallest instances that remained unsolved by different B\&C methods from the literature. For each number of customers and each number of vehicles, the highest filling coefficient $\bar{f}$ for which instances have been solved and the smallest filling coefficient leading to instances that could not be solved are reported. To increase readability, we only present what we consider to be the most important configurations per paper. We also indicate whether the customer demands were modeled as deterministic (D), Poisson RVs (P), or normal RVs (N) in each set of instances. Results from \cite{jabali2014partial} are omitted because they are slightly inferior to those of \cite{Ymro2023}.

Although the figure should be taken with a grain of salt, these results indicate that the DL-shaped method constitutes a significant advancement in the ability to solve instances with a high number of customers and vehicles, and a high filling coefficient. Although previous methods from the literature can solve some instances with a relatively high number of customers, they also fail to solve instances with as few as 30 or 40 customers.

\section{Conclusion}\label{sec:Conclusions}
This paper presents a new approach for solving a class of stochastic integer programs in which the first-stage solutions can be decomposed into disjoint components. The method is applied to the vehicle routing problem with stochastic demands under the detour-to-depot recourse policy. Our computational experiments show that it achieves state-of-the-art results on instances from the literature. Our approach is based on new optimality cuts and exploits new lower bounds on the recourse that are used by new lower bounding functionals. Computational experiments show that our lower bounds improve over the existing ones from the literature.

Regarding future lines of research, the method could be applied to other variants of the stochastic vehicle routing problem and other two-stage stochastic integer programs. Since the efficiency of this method depends on the availability of tight bounds on disjoint components of the recourse, this work may motivate the derivation of such bounds for other problems of interest. Finally, as our method requires the recourse function to be monotonic, it would be interesting to conduct a systematic study of this property for other two-stage stochastic programs.

\section*{Acknowledgments}
\label{sec:acknowledgments}
We thank Digital Research Alliance of Canada for providing high-performance computing facilities. Financial support for this work was provided by the Canadian Natural Sciences and Engineering Research Council (NSERC) under Grant 2021-04037. This support is gratefully acknowledged.

\bibliography{references}
\bibliographystyle{apalike} 

\end{document}